\documentclass{amsart}
\usepackage{amsmath,amssymb}
\usepackage[dvips]{graphics}
\usepackage[all]{xy}
\newtheorem{Theorem}{Theorem}[section]
\newtheorem{Lemma}[Theorem]{Lemma}
\newtheorem{Proposition}[Theorem]{Proposition}
\newtheorem{Corollary}[Theorem]{Corollary}

\newtheorem{Remark}[Theorem]{Remark}

\makeatletter
\@addtoreset{figure}{section}
\def\@thmcountersep{-}
\makeatother

\numberwithin{equation}{section}

\begin{document}

\title{Multiplicity of a space over another space}

\author{Kouki Taniyama}
\address{Department of Mathematics, School of Education, Waseda University, Nishi-Waseda 1-6-1, Shinjuku-ku, Tokyo, 169-8050, Japan}
\email{taniyama@waseda.jp}
\thanks{The author was partially supported by Grant-in-Aid for Scientific Research (C) (No. 21540000), Japan Society for the Promotion of Science.}


\date{}

\dedicatory{}


\begin{abstract}
We define a concept which we call multiplicity. First, multiplicity of a morphism is defined. Then the multiplicity of an object over another object is defined to be the minimum of the multiplicities of all morphisms from one to another. Based on this multiplicity, we define a pseudo distance on the class of objects. We define and study several multiplicities in the category of topological spaces and continuous maps, the category of groups and homomorphisms, the category of finitely generated $R$-modules and $R$-linear maps over a principal ideal domain $R$, and the neighbourhood category of oriented knots in the 3-sphere. 
\end{abstract}

\maketitle

\section{Introduction} 

Let ${\mathcal C}$ be a category with objects $X,Y,\cdots$. We denote the set of morphisms from $X$ to $Y$ by ${\rm Hom}(X,Y)$. By $f:X\to Y$ we mean $f\in{\rm Hom}(X,Y)$. The composition of $f:X\to Y$ and $g:Y\to Z$ is denoted by $g\circ f:X\to Z$. The identity morphism on $X$ is denoted by ${\rm id}_X:X\to X$. Note that by the definition of category the following (1), (2) and (3) hold.

\noindent
(1) For any $f:X\to Y$, $g:Y\to Z$ and $h:Z\to W$, $h\circ(g\circ f)=(h\circ g)\circ f$.

\noindent
(2) For any $f:X\to Y$, $f\circ {\rm id}_X=f$ and for any $g:Y\to X$, ${\rm id}_X\circ g=g$.

\noindent
(3) ${\rm Hom}(X,Y)$ and ${\rm Hom}(Z,W)$ are disjoint unless $X=Z$ and $Y=W$.

Let ${\mathbb R}$ be the set of all real numbers and ${\mathbb N}$ the set of all natural numbers. 
For a real number $a$ we denote the set of all real numbers greater than or equal to $a$ by ${\mathbb R}_{\geq a}$. Similarly we denote the set of all real numbers greater than $a$ by ${\mathbb R}_{> a}$. 
Let $\infty$ be an element that is not in ${\mathbb R}$. 
We extend the order, addition and multiplication of ${\mathbb R}$ to ${\mathbb R}\cup\{\infty\}$ in the usual way. Namely, for any real number $r$ it holds $r\leq\infty$, $\infty\leq\infty$, $r+\infty=\infty+r=\infty+\infty=\infty$ and for any positive real number $r$ it holds $r\cdot\infty=\infty\cdot r=\infty\cdot\infty=\infty$. 

Suppose that for each morphism $f:X\to Y$, an element $m(f)$ of ${\mathbb R}_{\geq 1}\cup\{\infty\}$ is assigned such that the following (1) and (2) hold.

\noindent
(1) $m({\rm id}_X)=1$ for any object $X$.

\noindent
(2) For any $f:X\to Y$ and $g:Y\to Z$, $m(g\circ f)\leq m(f)m(g)$.

\noindent
Then we say that $m$ is a {\it multiplicity} on the category ${\mathcal C}$. 
Let $m(X:Y)$ be the infimum of $m(f)$ where $f$ varies over all elements of ${\rm Hom}(X,Y)$. If there are no morphisms from $X$ to $Y$ then we define $m(X:Y)=\infty$. We call $m(X:Y)$ the {\it multiplicity of $X$ over $Y$}.

\vskip 3mm

\begin{Proposition}\label{Proposition 1.1}
{\rm(1)} For any objects $X$ and $Y$, $m(X:Y)\geq1$.

\noindent
{\rm(2)} For any object $X$, $m(X:X)=1$.

\noindent
{\rm(3)} For any objects $X$, $Y$ and $Z$,
\[
m(X:Z)\leq m(X:Y)m(Y:Z).
\]
\end{Proposition}

\vskip 3mm

\noindent{\bf Proof.} (1) By definition we have $m(X:Y)\geq1$. 

\noindent
(2) We have $m(X:X)=m({\rm id}_X)=1$. 

\noindent
(3) For any $\varepsilon>0$ there exists $f:X\to Y$ with 
\[
m(f)-m(X:Y)<{\rm min}\{\frac{\varepsilon}{3m(Y:Z)},\frac{\sqrt{3\varepsilon}}{3}\} 
\]
and $g:Y\to Z$ with 
\[
m(g)-m(Y:Z)<{\rm min}\{\frac{\varepsilon}{3m(X:Y)},\frac{\sqrt{3\varepsilon}}{3}\}.
\]
Then we have
\[
m(f)m(g)-m(X:Y)m(Y:Z)
\]
\[
=m(f)m(g)-m(f)m(Y:Z)+m(f)m(Y:Z)-m(X:Y)m(Y:Z)
\]
\[
=m(f)(m(g)-m(Y:Z))+(m(f)-m(X:Y))m(Y:Z)
\]
\[
<(m(X:Y)+\frac{\sqrt{3\varepsilon}}{3})(m(g)-m(Y:Z))+\frac{\varepsilon}{3m(Y:Z)}m(Y:Z)
\]
\[
<m(X:Y)\frac{\varepsilon}{3m(X:Y)}+\frac{\sqrt{3\varepsilon}}{3}\frac{\sqrt{3\varepsilon}}{3}+\frac{\varepsilon}{3m(Y:Z)}m(Y:Z)
\]
\[
=\varepsilon.
\]
Since $m(X:Z)\leq m(g\circ f)\leq m(f)m(g)$ we have $m(X:Z)-m(X:Y)m(Y:Z)<\varepsilon$. Since $\varepsilon$ is any positive number we have the conclusion. $\Box$

\vskip 3mm

\vskip 3mm

\begin{Remark}\label{discrete}
{\rm All examples of $m$ in this paper take their values in a proper subset ${\mathbb N}\cup\{\infty\}$ of ${\mathbb R}_{\geq 1}\cup\{\infty\}$ or a proper subset $\{e^n|n\in\{0\}\cup{\mathbb N}\}$ of ${\mathbb R}_{\geq 1}\cup\{\infty\}$ where $e$ is the base of natural logarithm.  
Then $m(X:Y)$ is simply defined as the minimum, not as the infimum, and the proof of Proposition \ref{Proposition 1.1} (3) becomes much simpler. However we define the range of $m$ as above for potential future use.}
\end{Remark}

\vskip 3mm

Let $d_m(X,Y)$ be an element of ${\mathbb R}_{\geq 0}\cup\{\infty\}$ defined by the following.
\[
d_m(X,Y)={\rm log}_e(m(X:Y)m(Y:X)).
\]
Here we define ${\rm log}_e(\infty)=\infty$ as usual. We call $d_m(X,Y)$ the {\it multiplicity distance} of $X$ and $Y$.

We say that $X$ has {\it finite multiplicity property} over $Y$ with respect to $m$ if $m(X:Y)\neq\infty$, namely if there is a morphism $f:X\to Y$ with $m(f)<\infty$. We say that $m$ has {\it finite multiplicity property} if any object $X$ has finite multiplicity over any object $Y$.

\vskip 3mm

\begin{Proposition}\label{pseudo distance}
Let $m$ be a multiplicity on a category ${\mathcal C}$ that has finite multiplicity property. Then $d_m$ is a pseudo distance on the class of objects of ${\mathcal C}$. Namely the following {\rm(D1$'$)}, {\rm(D2)} and {\rm(D3)} hold for any objects $X$, $Y$ and $Z$ of ${\mathcal C}$.

\noindent
{\rm(D1$'$)} $d_m(X,Y)\geq0$, $d_m(X,X)=0$,

\noindent
{\rm(D2)} $d_m(X,Y)=d_m(Y,X)$,

\noindent
{\rm(D3)} $d_m(X,Z)\leq d_m(X,Y)+d_m(Y,Z)$.

\end{Proposition}

\vskip 3mm

\noindent{\bf Proof.} (D1$'$) By Proposition \ref{Proposition 1.1} (1) we have $m(X:Y)m(Y:X)\geq1$. Therefore $d_m(X,Y)={\rm log}_e(m(X:Y)m(Y:X))\geq0$. By Proposition \ref{Proposition 1.1} (2) we have $m(X:X)m(X:X)=1$. Therefore $d_m(X,X)={\rm log}_e(m(X:X)m(X:X))=0$. 

\noindent
(D2) By definition we have  $d_m(X,Y)=d_m(Y,X)$.

\noindent
(D3) By Proposition \ref{Proposition 1.1} (3) we have $m(X:Z)m(Z:X)\leq m(X:Y)m(Y:Z)m(Z:Y)m(Y:X)=m(X:Y)m(Y:X)m(Y:Z)m(Z:Y)$. Therefore  $d_m(X,Z)={\rm log}_e(m(X:Z)m(Z:X))\leq{\rm log}_e(m(X:Y)m(Y:X)m(Y:Z)m(Z:Y))={\rm log}_e(m(X:Y)m(Y:X))+{\rm log}_e(m(Y:Z)m(Z:Y))=d_m(X,Y)+d_m(Y,Z)$.
$\Box$

\vskip 3mm

The following is a typical example of multiplicity. Let ${\mathcal C}_{\rm set}$ be the category of non-empty sets and maps. Then a morphism $f:X\to Y$ is a map from a set $X$ to a set $Y$. We denote the cardinality of a set $A$ by $|A|$. For any infinite set $A$, the cardinality of $A$ is denoted by the same symbol $|A|=\infty$. 
Let $m_{\rm map}(f)$ be an element of ${\mathbb N}\cup\{\infty\}$ defined by
\[
m_{\rm map}(f)={\rm sup}\{|f^{-1}(y)|\mid y\in Y\}.
\]
We call $m_{\rm map}$ {\it map-multiplicity}. 

\vskip 3mm

\begin{Proposition}\label{map-multiplicity} A map-multiplicity $m_{\rm map}$ is a multiplicity. Namely the following (1) and (2) hold.

\noindent
(1) $m_{\rm map}({\rm id}_X)=1$ for any set $X$ of ${\mathcal C}_{\rm set}$.

\noindent
(2) For any maps $f:X\to Y$ and $g:Y\to Z$, $m_{\rm map}(g\circ f)\leq m_{\rm map}(f)m_{\rm map}(g)$.
\end{Proposition}

\vskip 3mm

\noindent{\bf Proof.} The proof of (1) is trivial. We will show (2). If $m_{\rm map}(f)=\infty$ or $m_{\rm map}(g)=\infty$ then the inequality holds. Suppose $m_{\rm map}(f)<\infty$ and $m_{\rm map}(g)<\infty$. Let $z$ be any element of $Z$. Since $|g^{-1}(z)|\leq m_{\rm map}(g)$ and $|f^{-1}(y)|\leq m_{\rm map}(f)$ for any $y\in Y$ we have $|(g\circ f)^{-1}(z)|=|\cup_{y\in g^{-1}(z)}f^{-1}(y)|\leq m_{\rm map}(f)m_{\rm map}(g)$. Thus we have $m_{\rm map}(g\circ f)\leq m_{\rm map}(f)m_{\rm map}(g)$.
$\Box$

\vskip 3mm

The following proposition shows that map-multiplicity is, in a sense, generalization of division. 
For a real number $a$ we denote the least integer no less than $a$ by $\lceil a\rceil$. 

\vskip 3mm

\begin{Proposition}\label{division} 
Let $X$ and $Y$ be finite sets. Then the following holds. 
\[
m_{\rm map}(X:Y)=\lceil|X|/|Y|\rceil.
\]
\end{Proposition}

\vskip 3mm

\noindent{\bf Proof.} 
For any map $f:X\to Y$ we have $m_{\rm map}(f)\geq\lceil|X|/|Y|\rceil$ by the pigeonhole principle, and clearly, there is a map $g:X\to Y$ with $m_{\rm map}(g)=\lceil|X|/|Y|\rceil$. Therefore we have the result. 
$\Box$

\vskip 3mm

Let ${\mathcal C}_1$ and ${\mathcal C}_2$ be categories and $m_2$ a multiplicity on ${\mathcal C}_2$. Let $F$ be a functor from ${\mathcal C}_1$ to ${\mathcal C}_2$. For each morphism $f$ of ${\mathcal C}_1$ we define $m_1(f)$ by $m_1(f)=m_2(F(f))$. We call $m_1$ the {\it pull-back multiplicity} of $m_2$ with respect to a functor $F$. 

\vskip 3mm

\begin{Proposition}\label{pull-back multiplicity} The pull-back multiplicity $m_1$ is a multiplicity. Namely the following (1) and (2) hold.

\noindent
(1) $m_1({\rm id}_X)=1$ for any object $X$ of ${\mathcal C}_1$.

\noindent
(2) For any morphisms $f:X\to Y$ and $g:Y\to Z$ of ${\mathcal C}_1$, $m_1(g\circ f)\leq m_1(f)m_1(g)$.
\end{Proposition}

\vskip 3mm

\noindent{\bf Proof.} (1) By definition we have $m_1({\rm id}_X)=m_2(F({\rm id}_X))$. Since $F$ is a functor we have $F({\rm id}_X)={\rm id}_{F(X)}$. Therefore $m_2(F({\rm id}_X))=m_2({\rm id}_{F(X)})=1$.

\noindent
(2) By definition we have $m_1(g\circ f)=m_2(F(g\circ f))$. Since $F$ is a functor we have $F(g\circ f)=F(g)\circ F(f)$. 
Therefore $m_2(F(g\circ f))=m_2(F(g)\circ F(f))\leq m_2(F(f))m_2(F(g))=m_1(f)m_1(g)$. 
$\Box$

\vskip 3mm

Suppose that a category has a functor to the category of non-empty sets and maps, then the pull-back multiplicity of the map-multiplicity is also called the map-multiplicity so long as no confusion occurs. Many categories, such as the category of topological spaces and continuous maps and the category of groups and homomorphisms, have a forgetful functor to the category of non-empty sets and maps. Thus the map-multiplicity is defined on such categories. We will study them in Section \ref{topology} and Section \ref{group}. 

Let $X$ and $Y$ be objects. Suppose that there exists an injective morphism from $X$ to $Y$, and there exists an injective morphism from $Y$ to $X$. Then we consider the problem whether or not $X$ and $Y$ are isomorphic. We consider this problem in several categories in the following sections. An analogous problem for surjective morphism is also considered in Section \ref{group} and Section \ref{module}.

\vskip 3mm

\begin{Remark}\label{historical remark}
{\rm In \cite{K-T} a natural number $b(f)$ is defined for an embedding $f:G\to{\mathbb S}^3$ of a finite graph $G$ into the $3$-sphere ${\mathbb S}^3$ as a generalization of the braid index of knots and links. Then $b(f)$ is estimated below by $m_{\rm map}(g)$ for any continuous map $g:G\to{\mathbb S}^1$ from $G$ to the unit circle ${\mathbb S}^1$. Therefore the author defined map-multiplicity of continuous maps and then generalized it to multiplicities in various categories. There are some preliminary announcements of this work. See \cite{Taniyama}, \cite{Taniyama2} and \cite{Taniyama3}. 
In \cite{B-F-K} S. Bogatyi, J. Fricke and E. Kudryavtseva independently defined the same number as our map-multiplicity of a continuous map, and studied it in somewhat different interest. In \cite{Gromov} M. Gromov independently defined the cardinality of a topological space over another topological space which is essentially the same as our map-multiplicity of a topological space over another topological space, and studied it in various aspects. In \cite{Karasev} R. Karasev also showed some results on multiplicity of continuous maps between manifolds.
}
\end{Remark}

\vskip 3mm

\section{Multiplicity of topological spaces}\label{topology} 

Let $X$ and $Y$ be topological spaces and $f:X\to Y$ a continuous map from $X$ to $Y$. Then the map-multiplicity of $f$ is defined by 
$m_{\rm map}(f)={\rm sup}\{|f^{-1}(y)|\mid y\in Y\}$ and the map-multiplicity $m_{\rm map}(X:Y)$ of $X$ over $Y$ is defined to be the infimum of $m_{\rm map}(f)$ where $f$ varies over all continuous maps from $X$ to $Y$. 

We say that topological spaces $X$ and $Y$ are {\it weakly homeomorphic} if there exists a continuous injection from $X$ to $Y$ and there exists a continuous injection from $Y$ to $X$. We say that a set of topological spaces ${\mathcal T}$ is {\it classed} if any two weakly homeomorphic topological spaces in $\mathcal T$ are homeomorphic. We say that a topological space $X$ is {\it self-closed} if every continuous injection from $X$ to $X$ is a homeomorphism.

\vskip 3mm

\begin{Proposition}\label{self-closed topological space}
A set of self-closed topological spaces is classed.
\end{Proposition}

\vskip 3mm

\noindent{\bf Proof.} Let $X$ and $Y$ be self-closed topological spaces. Suppose that $X$ and $Y$ are weakly homeomorphic. Then there are continuous injections $f:X\to Y$ and $g:Y\to X$. Then $g\circ f:X\to X$ is a continuous injection. Since $X$ is self-closed $g\circ f$ is a homeomorphism. Similarly $f\circ g:Y\to Y$ is a homeomorphism. Therefore both $f$ and $g$ are bijections. Since $f^{-1}:Y\to X$ is a composition of two continuous maps $g:Y\to X$ and $(g\circ f)^{-1}:X\to X$, $f^{-1}$ is also continuous. Therefore $f$ is a homeomorphism. $\Box$

\vskip 3mm

\begin{Theorem}\label{classed topological spaces}
Let $\mathcal{X}$ be a classed set of topological spaces. Let ${\mathcal C}_{\mathcal{X}}$ be a category whose objects are the elements of $\mathcal{X}$ and morphisms are the continuous maps between the elements of $\mathcal{X}$. Let $m_{\rm map}$ be the map-multiplicity on ${\mathcal C}_{\mathcal{X}}$ and $d_{m_{\rm map}}$ the map-multiplicity distance. Let $X$ and $Y$ be elements of $\mathcal{X}$. Then $X$ and $Y$ are homeomorphic if and only if $d_{m_{\rm map}}(X,Y)=0$. Thus, if $m_{\rm map}$ has finite multiplicity property on $\mathcal{X}$, then the pseudo distance $d_{m_{\rm map}}$ defines a distance on the set of homeomorphism classes of $\mathcal{X}$.
\end{Theorem}

\vskip 3mm

\noindent{\bf Proof.} Suppose that $X$ and $Y$ are homeomorphic. Then there is a homeomorphism $f:X\to Y$. Since $m_{\rm map}(f)=m_{\rm map}(f^{-1})=1$ we have $m_{\rm map}(X:Y)=m_{\rm map}(Y:X)=1$. Therefore $d_{m_{\rm map}}(X,Y)=0$. Suppose $d_{m_{\rm map}}(X,Y)=0$. Then there are continuous injections $f:X\to Y$ and $g:Y\to X$. Therefore $X$ and $Y$ are weakly homeomorphic. Since $\mathcal{X}$ is classed, $X$ and $Y$ are homeomorphic. $\Box$

\vskip 3mm

In the following we consider simplicial complexes. All simplicial complexes in this paper are locally finite. 
Note that a simplicial complex is formally a set of simplexes. However we do not distinguish a simplicial complex and the union of its simplices so long as no confusion occurs. 

\vskip 3mm

\begin{Lemma}\label{graph-subspace}
Let $G$ be a finite one-dimensional simplicial complex. Let $X$ be a subspace of $G$ that is homeomorphic to a finite one-dimensional complex. 
Then the following (1), (2) and (3) hold. 

\noindent
(1) For any $1$-simplex $e$ of $G$, $e\cap X$ is a disjoint union of finitely many points and closed intervals. 

\noindent
(2) There is a subdivision $G'$ of $G$ such that $X$ is a sub-complex of $G'$.

\noindent
(3) The first Betti number $\beta_1(X)$ of $X$ is less than or equal to the first Betti number $\beta_1(G)$ of $G$.
\end{Lemma}

\vskip 3mm

\noindent{\bf Proof.} 
(1) Since $X$ is compact and $G$ is Hausdorff, $X$ is closed in $G$. Also $e$ is closed in $G$. Thus $e\cap X$ is a closed subset of a compact set $G$, hence compact. Each connected component of $e\cap X$ is a compact connected subset of $e$. Therefore it is a point or a closed interval. Suppose that the number of connected components of $e\cap X$ is not finite. Since $X$ is compact the number of connected components of $X$ is finite. Therefore there are connected components $A_1,A_2$ and $A_3$ of $e\cap X$ that are contained in a connected component $B$ of $X$. We may suppose without loss of generality that $A_1,A_2$ and $A_3$ are arranged in this order in $e$. Then there are points $x$ and $y$ of $e\setminus X$ such that $x$ is between $A_1$ and $A_2$ and $y$ is between $A_2$ and $A_3$. But then $A_2$ is separated from $A_1$ and $A_3$ in $X$ and they cannot be contained in the same connected component of $X$. This is a contradiction. Therefore the number of connected components of $e\cap X$ is finite. Thus $e\cap X$ is a disjoint union of finitely many points and closed intervals. 

\noindent
(2) By (1) we immediately have the result.

\noindent
(3) By (2), $X$ is a sub-complex of a finite one-dimensional simplicial complex $G'$. Then the inequality $\beta_1(X)\leq\beta_1(G')$ is a well-known fact in combinatorial topology. Since $\beta_1(G')=\beta_1(G)$ we have the result.
$\Box$

\vskip 3mm

\begin{Proposition}\label{self-closed spaces}

\noindent
(1) A compact Hausdorff space $X$ is self-closed if and only if no proper subspace of $X$ is homeomorphic to $X$.

\noindent
(2) A closed manifold is self-closed.

\noindent
(3) Let $G$ be a finite one-dimensional simplicial complex. Suppose that no $0$-simplex is a face of exactly one $1$-simplex. Then $G$ is self-closed.
\end{Proposition}

\vskip 3mm

\noindent{\bf Proof.} (1) Suppose that a proper subset $Y$ of $X$ is homeomorphic to $X$. Let $f:X\to Y$ be a homeomorphism and $i:Y\to X$ the inclusion map. Then the composition $i\circ f:X\to X$ is a continuous injection. Since $i\circ f$ is not surjective it is not a homeomorphism. Therefore $X$ is not self-closed. Suppose that $X$ is not self-closed. Then there is a continuous injection $f:X\to X$ that is not a homeomorphism. Since the source $X$ is compact and the target $X$ is Hausdorff, $f$ is a closed map. Therefore $X$ and $f(X)$ are homeomorphic. Since $f$ is not a homeomorphism, $f(X)$ is a proper subspace of $X$.

\noindent
(2) Let $M$ be a closed $n$-dimensional manifold. Since $M$ is compact Hausdorff it is sufficient to show that no proper subspace of $M$ is homeomorphic to $M$. Let $M_1,\cdots,M_k$ be the connected components of $M$. Let $N$ be a subspace of $M$ that is homeomorphic to $M$. Set $N_i=N\cap M_i$. Since $N_i$ is compact and $M_i$ is Hausdorff, $N_i$ is a closed subset of $M_i$. 
We will show that $N_i$ is an open subset of $M_i$ by using the invariance of domain theorem. 
Let $y$ be a point in $N_i$. We will show that there is an open set $O$ of $M_i$ with $y\in O\subset N_i$. Let $U$ be an open neighbourhood of $y$ in $M_i$ that is homeomorphic to an open set $A$ of ${\mathbb R}^n$. Let $V$ be an open neighbourhood of $y$ in $N_i$ that is homeomorphic to an open set $B$ of ${\mathbb R}^n$. Let $f:A\to U$ and $g:B\to V$ be homeomorphisms. Since $U\cap V$ is an open subset of $V$, $g^{-1}(U\cap V)$ is an open set of ${\mathbb R}^n$. By the invariance of domain theorem we see that $f^{-1}(U\cap V)$ is an open set of ${\mathbb R}^n$. Therefore $f^{-1}(U\cap V)$ is an open set of $A$. Since $f$ is a homeomorphism $U\cap V$ is an open set of $U$. Therefore $U\cap V$ is an open set of $M_i$. Thus $O=U\cap V$ is the desired set. 
Since $M_i$ is connected we have $N_i=\emptyset$ or $N_i=M_i$. Since $N$ has the same number of connected components as $M$ we have $N_i=M_i$ for every $i$. Therefore $N=M$. 

\noindent

\noindent
(3) Since $G$ is compact Hausdorff, it is sufficient to show that no proper subspace of $G$ is homeomorphic to $G$. Let $H$ be a subspace of $G$ homeomorphic to $G$. For each $1$-simplex $e$ of $G$, $e\cap H$ is a disjoint union of finitely many points and closed intervals by Lemma \ref{graph-subspace}. If there is a $1$-simplex $e$ of $G$ with $e\cap H\neq e$, then by the assumption of $G$ we have $\beta_1(H)<\beta_1(G)$ or $\beta_0(G\setminus{\rm int}(e))=\beta_0(G)+1$ where $\beta_0(X)$ is the number of connected component of $X$. Then it is easy to see that $H$ cannot be homeomorphic to $G$.  Thus $e\cap H=e$ for each $1$-simplex $e$ of $G$. Since $\beta_0(H)=\beta_0(G)$, each $0$-simplex of $G$ is contained in $H$. Thus we have $H=G$ as desired. $\Box$ 

\vskip 3mm

Let $G$ be a one-dimensional simplicial complex. Let $x$ be a point of $G$. The {\it degree} of $x$ in $G$, denoted by ${\rm deg}(x)={\rm deg}(x,G)$, is defined as follows. If $x$ is a $0$-simplex of $G$, then ${\rm deg}(x)$ is defined to be the number of $1$-simplexes of $G$ that contain $x$. If $x$ is an interior point of a $1$-simplex of $G$, then ${\rm deg}(x)=2$. Note that ${\rm deg}(x,G')={\rm deg}(x,G)$ for any subdivision $G'$ of $G$. Let ${\rm br}(G)$ be the set of all points $x$ of $G$ with ${\rm deg}(x)={\rm deg}(x,G)\geq3$. 

\vskip 3mm

\begin{Lemma}\label{degree}
Let $f:X\to Y$ be a homeomorphism from a one-dimensional simplicial complex $X$ to another one-dimensional simplicial complex $Y$. Let $x$ be a point of $X$. Then ${\rm deg}(x,X)={\rm deg}(f(x),Y)$. Therefore $f({\rm br}(X))={\rm br}(Y)$.
\end{Lemma}

\vskip 3mm

\noindent{\bf Proof.} 
Suppose that ${\rm deg}(x,X)<{\rm deg}(f(x),Y)$. Let $M$ be a sufficiently small closed connected neighbourhood of $f(x)$ in $Y$. Then the number of connected components of $M\setminus\{f(x)\}$ is ${\rm deg}(f(x),Y)$. Let $N$ be a sufficiently small closed connected neighbourhood of $x$ in $X$ such that $f(N)\subset M$. Then the number of connected components of $N\setminus\{x\}$ is ${\rm deg}(x,X)$. 
Since $f$ is a homeomorphism $f(N)$ is also a neighbourhood of $f(x)$ in $Y$. Each connected component of $N\setminus\{x\}$ is mapped into a connected component of $M\setminus\{f(x)\}$ under $f$. Therefore there is a component of $M\setminus\{f(x)\}$ that is disjoint from $f(N)$. Then $f(N)$ cannot be a neighbourhood of $f(x)$ in $Y$. This is a contradiction. Thus we have ${\rm deg}(x,X)\geq{\rm deg}(f(x),Y)$. Similarly, by considering $f^{-1}$ we have ${\rm deg}(x,X)\leq{\rm deg}(f(x),Y)$.
$\Box$

\vskip 3mm

\begin{Proposition}\label{classed set of graphs}
Let ${\mathcal S}$ be a set of finite one-dimensional simplicial complexes. Suppose that each element of ${\mathcal S}$ has no isolated points and no connected components each of which is homeomorphic to a closed interval. Then ${\mathcal S}$ is classed.
\end{Proposition}

\vskip 3mm

\noindent{\bf Proof.} Let $G_1$ and $G_2$ be elements of ${\mathcal S}$. 
Suppose that there are continuous injections $f_1:G_1\to G_2$ and $f_2:G_2\to G_1$. 
Let $H_i$ be the maximal sub-complex of $G_i$ such that each $0$-simplex of $H_i$ is a face of at least two $1$-simplexes of $H_i$ for $i=1,2$. Namely $H_i$ is the union of all subspaces of $G_i$ each of which is homeomorphic to a circle. If $\beta_1(G_i)=0$, then $H_i$ is an empty set. Then $\beta_1(H_i)=\beta_1(G_i)$ for $i=1,2$. We may suppose without loss of generality that ${\rm deg}(x,G_i)\neq2$ for each $0$-simplex $x$ contained in $G_i\setminus H_i$. 
Since $G_1$ and $G_2$ are compact Hausdorff, $f_1$ and $f_2$ are embeddings and $G_i$ and $f_i(G_i)$ are homeomorphic for $i=1,2$. 
By a standard argument of combinatorial topology we have $\beta_1(f_1(G_1))=\beta_1(f_1(G_1)\cap H_2)$. 
Suppose that $H_2$ is not contained in $f_1(G_1)$. Then $f_1(G_1)\cap H_2$ is a proper subspace of $H_2$. Then we have $\beta_1(f_1(G_1)\cap H_2)<\beta_1(H_2)$ by a standard argument of combinatorial topology. Thus $\beta_1(G_1)=\beta_1(f_1(G_1))=\beta_1(f_1(G_1)\cap H_2)$ is less than $\beta_1(H_2)=\beta_1(G_2)$. 
On the other hand $\beta_1(G_2)=\beta_1(f_2(G_2))\leq\beta_1(G_1)$ by Lemma \ref{graph-subspace}. This is a contradiction. Therefore $H_2$ is contained in $f_1(G_1)$. 
Since $f_1(G_1)$ is a sub-complex of a subdivision of $G_2$, ${\rm br}(f_1(G_1))$ is a subset of ${\rm br}(G_2)$. Since ${\rm br}(f_1(G_1))=f_1({\rm br}(G_1))$ we have $|{\rm br}(G_1)|\leq|{\rm br}(G_2)|$. Similarly we have $|{\rm br}(G_2)|\leq|{\rm br}(G_1)|$. Therefore  ${\rm br}(f_1(G_1))={\rm br}(G_2)$. Similarly we have ${\rm br}(f_2(G_2))={\rm br}(G_1)$. 
Let $e$ be a $1$-simplex of $G_2$ that is not contained in $H_2$. Let $x$ and $y$ be the $0$-simplexes of $G_2$ contained in $e$. Suppose that both ${\rm deg}(x,G_2)$ and ${\rm deg}(y,G_2)$ are greater than $2$. We will show that $e$ is contained in $f_1(G_1)$. Suppose that $e$ is not contained in $f_1(G_1)$. 
By the condition ${\rm br}(f_1(G_1))={\rm br}(G_2)$ there is an interior point $z$ of $e$ that is not contained in $f_1(G_1)$. Note that $\beta_0(G_2\setminus\{z\})=\beta_0(G_2)+1$. Each connected component of $G_2\setminus\{z\}$ contains at least one point in ${\rm br}(G_2)$. Therefore $\beta_0(f_1(G_1))\geq\beta_0(G_2\setminus\{z\})$. Thus we have $\beta_0(G_1)=\beta_0(f_1(G_1))$ is greater than $\beta_0(G_2)$. 
Since each connected component of $G_1$ contains at least one point in ${\rm br}(G_1)$ and ${\rm br}(f_2(G_2))={\rm br}(G_1)$, we have $\beta_0(G_2)=\beta_0(f_2(G_2))\geq\beta_0(G_1)$. This is a contradiction. Thus we have $e\subset f_1(G_1)$. Let $e$ be a $1$-simplex of $G_2$ such that one of the $0$-simplexes of $G_2$ contained in $e$, say $x$, has degree $1$ in $G_2$. Note that by the assumption on $G_2$ the other $0$-simplex of $G_2$ contained in $e$, say $y$, has degree greater than $2$. Therefore $y\in{\rm br}(G_2)$. Therefore $e\cap f_1(G_1)$ is a neighbourhood of $y$ in $e$. By Lemma \ref{graph-subspace} and by the assumption on $G_1$, $e\cap f_1(G_1)$ is homeomorphic to a closed interval. Therefore $f_1(G_1)$ and $G_2$ are homeomorphic. 
$\Box$

\vskip 3mm

\begin{Proposition}\label{simplicial-complexes}
Let $n$ be a non-negative integer. Let $X$ and $Y$ be $n$-dimensional finite simplicial complexes. Then $X$ has finite multiplicity property over $Y$ with respect to map-multiplicity.
\end{Proposition}

\vskip 3mm

\noindent{\bf Proof.} It is sufficient to show that there is a continuous map from $f:X\to Y$ such that $f^{-1}(y)$ is a finite set for every $y\in Y$. Let $v_1,\cdots,v_k$ be the $0$-simplexes of $X$. Let $s$ be an $n$-simplex of $Y$. Let $x_1,\cdots,x_k$ be points in $s$ that are in general position. Namely any $l$ points of them with $l\leq n+1$ are not contained in any $l-2$-dimensional affine subspace. Let $f:X\to Y$ be a map such that $f(v_i)=x_i$ for each $i$ and the restriction map of $f$ to each simplex of $X$ is an affine map. Then $f$ has the desired property.
$\Box$

\vskip 3mm

The following is an example of map-multiplicity of a topological space over another topological space. See \cite{Gromov} and \cite{Karasev} for other examples. 

\begin{Proposition}\label{complete-graph}
Let $n$ be a natural number and $K_{2n+1}$ the $1$-skeleton of a $2n$-simplex. Let ${\mathbb S}^1$ be the unit circle. Then the following holds. 
\[
m_{\rm map}(K_{2n+1}:{\mathbb S}^1)=\frac{n(n+1)}{2}.
\]
\end{Proposition}

\vskip 3mm

\noindent{\bf Proof.} Let $f:K_{2n+1}\to {\mathbb S}^1$ be the map constructed as follows. The image of the $0$-simplexes of $K_{2n+1}$ under $f$ forms the set of the vertices of a regular $(2n+1)$-gon inscribing to ${\mathbb S}^1$. For each $1$-simplex $e$ of $K_{2n+1}$, $f$ maps $e$ homeomorphically onto an arc in ${\mathbb S}^1$ that is shorter than another arc in ${\mathbb S}^1$ with the same end points. Then it is easy to check that $\displaystyle{m_{\rm map}(f)=\frac{n(n+1)}{2}}$. Thus we have $\displaystyle{m_{\rm map}(K_{2n+1}:{\mathbb S}^1)\leq\frac{n(n+1)}{2}}$. We will show $\displaystyle{m_{\rm map}(g)\geq\frac{n(n+1)}{2}}$ for any continuous map $g:K_{2n+1}\to {\mathbb S}^1$. If $g$ maps a $1$-simplex to a point, then $m_{\rm map}(g)=\infty$. Thus we may suppose that $g$ maps no $1$-simplex to a point. We will deform $g$ step by step without increasing $m_{\rm map}(g)$ as follows. 
Let $e$ be a $1$-simplex of $K_{2n+1}$ and $u$ and $v$ the $0$-simplexes contained in $e$. 
First suppose $g(u)\neq g(v)$. Since $g(e)$ is a connected subset of ${\mathbb S}^1$, $g(e)$ contains an arc, say $\alpha$, of ${\mathbb S}^1$ joining $g(u)$ and $g(v)$. Then we modify $g$ on $e$, still denoted by $g$, such that $g$ maps $e$ homeomorphically onto $\alpha$. 
Next suppose that $g(u)=g(v)$. Since $g$ does not map $e$ to a point, there is an interior point, say $x$, of $e$ such that $g(u)\neq g(x)$. Let $e_1$ be the line segment in $e$ with end points $u$ and $x$, and $e_2$ the line segment in $e$ with end points $x$ and $v$. 
Since $g(e_i)$ is a connected subset of ${\mathbb S}^1$, $g(e_i)$ contains an arc, say $\beta_i$, of ${\mathbb S}^1$ joining $g(u)=g(v)$ and $g(x)$ for $i=1,2$. Then we modify $g$ on $e$ such that $g$ maps $e_i$ homeomorphically onto $\beta_i$ for $i=1,2$. We perform this modification for each $1$-simplex of $K_{2n+1}$. Thus there is a subdivision, say $K_{2n+1}'$ of $K_{2n+1}$ such that $g$ maps each $1$-simplex of $K_{2n+1}'$ homeomorphically onto an arc in ${\mathbb S}^1$. Then by a slight modification near the $0$-simplexes of $K_{2n+1}'$ we can modify $g$ without increasing $m_{\rm map}(g)$ such that no two $0$-simplexes of $K_{2n+1}'$ are mapped to the same point under $g$. Then we can deform $g$ such that $g$ maps each $1$-simplex of $K_{2n+1}$ homeomorphically onto an arc in ${\mathbb S}^1$. Thus $g$ maps the $0$-simplexes of $K_{2n+1}$ to $2n+1$ points in ${\mathbb S}^1$ and each $1$-simplex of $K_{2n+1}$ homeomorphically onto an arc in ${\mathbb S}^1$. Let $V$ be the set of $0$-simplexes of $K_{2n+1}$. Let $e$ be a $1$-simplex of $K_{2n+1}$. Let $d(e,g)=|g(e)\cap g(V)|-1$. Let $v$ be a $0$-simplex of $K_{2n+1}$. Then the sum of $d(e,g)$ over all $1$-simplexes $e$ of $K_{2n+1}$ containing $v$ is greater than or equal to $2(1+2+\cdots+n)$. Therefore the sum $S$ of $d(e,g)$ over all $1$-simplexes $e$ of $K_{2n+1}$ is greater than or equal to $\displaystyle{\frac{2(1+2+\cdots+n)(2n+1)}{2}}=\frac{n(n+1)(2n+1)}{2}$. 
Let $y_1,\cdots,y_{2n+1}$ be points in ${\mathbb S}^1\setminus g(V)$ such that each connected component of ${\mathbb S}^1\setminus g(V)$ contains exactly one of $y_1,\cdots,y_{2n+1}$. Then $d(e,g)=|g(e)\cap\{y_1,\cdots,y_{2n+1}\}|$. 
Note that $|g^{-1}(y_i)|$ is equal to the number of $1$-simplexes $e$ of $K_{2n+1}$ with $g(e)\ni y_i$. 
Therefore the sum $|g^{-1}(y_1)|+\cdots+|g^{-1}(y_{2n+1})|$ is equal to $S$. 
Then by the pigeonhole principle there is a point $y_i$ such that $\displaystyle{|g^{-1}(y_i)|\geq\frac{n(n+1)(2n+1)}{2(2n+1)}=\frac{n(n+1)}{2}}$. Thus we have $\displaystyle{m_{\rm map}(g)\geq\frac{n(n+1)}{2}}$. 
$\Box$

\vskip 3mm

\section{Multiplicity of groups}\label{group} 

\vskip 3mm

Let $G$ and $H$ be groups and $f:G\to H$ a homomorphism from $G$ to $H$. 
Then the map-multiplicity of $f$ is defined by $m_{\rm map}(f)={\rm sup}\{|f^{-1}(y)|\mid y\in H\}$. 
Note that $m_{\rm map}(f)=|{\rm Ker}(f)|$ where $|{\rm Ker}(f)|$ is the order of the kernel ${\rm Ker}(f)$ of $f$. 
Therefore we denote $m_{\rm map}(f)$ by $m_{\rm Ker}(f)$ and call it {\it kernel-multiplicity} of $f$ in this section. 
Then the kernel-multiplicity $m_{\rm Ker}(G:H)$ of $G$ over $H$ is defined to be the infimum of $m_{\rm Ker}(f)$ where $f$ varies over all homomorphisms from $G$ to $H$. 

We say that groups $G$ and $H$ are {\it weakly isomorphic} if there exists an injective homomorphism from $G$ to $H$ and there exists an injective homomorphism from $H$ to $G$. We say that a set of groups ${\mathcal G}$ is {\it classed} if any two weakly isomorphic groups in $\mathcal G$ are isomorphic. A group $G$ is said to be {\it co-Hopfian} if every injective homomorphism from $G$ to $G$ is an isomorphism. 

\vskip 3mm

\begin{Proposition}\label{co-Hopfian}
A set of co-Hopfian groups is classed.
\end{Proposition}

\vskip 3mm

\noindent{\bf Proof.} 
Suppose that $G$ and $H$ are mutually weakly isomorphic co-Hopfian groups. Let $f:G\to H$ and $g:H\to G$ be injective homomorphisms. Then $g\circ f:G\to G$ is an injective homomorphism. Since $G$ is co-Hopfian, $g\circ f$ is an isomorphism. Hence $g$ is surjective and therefore $g$ is an isomorphism from $H$ to $G$. 
$\Box$

\vskip 3mm

\begin{Theorem}\label{classed groups}
Let ${\mathcal G}$ be a classed set of groups. Let ${\mathcal C}_{\mathcal G}$ be a category whose objects are the elements of ${\mathcal G}$ and morphisms are the homomorphisms between the elements of ${\mathcal G}$. Let $m_{\rm Ker}$ be the kernel-multiplicity on ${\mathcal C}_{\mathcal G}$ and $d_{m_{\rm Ker}}$ the kernel-multiplicity distance. Let $G$ and $H$ be elements of ${\mathcal G}$. Then $G$ and $H$ are isomorphic if and only if $d_{m_{\rm Ker}}(G,H)=0$. Thus, if $m_{\rm Ker}$ has finite multiplicity property on ${\mathcal G}$, then the pseudo distance $d_{m_{\rm Ker}}$ defines a distance on the set of isomorphism classes of ${\mathcal G}$.
\end{Theorem}

\vskip 3mm

The proof of Theorem \ref{classed groups} is entirely analogous to that of Theorem \ref{classed topological spaces} and we omit it. Note that a finite group is co-Hopfian, and kernel-multiplicity has finite multiplicity property on the set of all finite groups. 
Thus $d_{m_{\rm Ker}}$ defines a distance on the set of isomorphism classes of all finite groups. 

We now define another multiplicity for group homomorphisms. Let $G$ and $H$ be groups and $f:G\to H$ a homomorphism from $G$ to $H$. Then the cokernel of $f$ is the quotient set ${\rm Coker}(f)=H/f(G)$. Then we call $m_{\rm Coker}(f)=|{\rm Coker}(f)|$ the {\it cokernel-multiplicity} of $f$.

\vskip 3mm

\begin{Proposition}\label{cokernel-multiplicity} A cokernel-multiplicity $m_{\rm Coker}$ is a multiplicity. Namely the following (1) and (2) hold.

\noindent
(1) $m_{\rm Coker}({\rm id}_G)=1$ for any group $G$.

\noindent
(2) For any group homomorphisms $f:G\to H$ and $g:H\to K$, $m_{\rm Coker}(g\circ f)\leq m_{\rm Coker}(f)m_{\rm Coker}(g)$.
\end{Proposition}

\vskip 3mm

\noindent{\bf Proof.} (1) Since ${\rm Coker}({\rm id}_G)=G/G$ is a trivial group, $m_{\rm Coker}({\rm id}_G)=|G/G|=1$.

\noindent
(2) Since ${\rm Coker}(g\circ f)=K/g\circ f(G)=K/g(f(G))$, $m_{\rm Coker}(g\circ f)=|{\rm Coker}(g\circ f)|=|K/g(f(G))|=|K/g(H)|\cdot|g(H)/g(f(G))|$. 
Since $|K/g(H)|=m_{\rm Coker}(g)$ it is sufficient to show $|g(H)/g(f(G))|\leq m_{\rm Coker}(f)=|H/f(G)|$. There is a surjection from $H/f(G)$ to $g(H)/g(f(G))$ induced by $g$. Therefore $|g(H)/g(f(G))|\leq|H/f(G)|$. 
$\Box$

\vskip 3mm

The following arguments show an analogy of cokernel-multiplicity to kernel-multiplicity. The proofs are entirely analogous and we omit them. 
We say that groups $G$ and $H$ are {\it co-weakly isomorphic} if there exists a surjective homomorphism from $G$ to $H$ and there exists a surjective homomorphism from $H$ to $G$. We say that a set of groups ${\mathcal G}$ is {\it co-classed} if any two co-weakly isomorphic groups in $\mathcal G$ are isomorphic. A group $G$ is said to be {\it Hopfian} if every surjective homomorphism from $G$ to $G$ is an isomorphism. 

\vskip 3mm

\begin{Proposition}\label{Hopfian}
A set of Hopfian groups is co-classed.
\end{Proposition}

\vskip 3mm

\begin{Theorem}\label{co-classed groups}
Let ${\mathcal G}$ be a co-classed set of groups. Let ${\mathcal C}_{\mathcal G}$ be a category whose objects are the elements of ${\mathcal G}$ and morphisms are the homomorphisms between the elements of ${\mathcal G}$. Let $m_{\rm Coker}$ be the cokernel-multiplicity on ${\mathcal C}_{\mathcal G}$ and $d_{m_{\rm Coker}}$ the cokernel-multiplicity distance. Let $G$ and $H$ be elements of ${\mathcal G}$. Then $G$ and $H$ are isomorphic if and only if $d_{m_{\rm Coker}}(G,H)=0$. Thus, if $m_{\rm Coker}$ has finite multiplicity property on ${\mathcal G}$, then the pseudo distance $d_{m_{\rm Coker}}$ defines a distance on the set of isomorphism classes of ${\mathcal G}$.
\end{Theorem}

\vskip 3mm

The following example shows a difference between kernel-multiplicity and cokernel-multiplicity. Let $F(n)$ be a free group of rank $n$. 
Then there exists an injective homomorphism from $F(2)$ to $F(3)$, and a surjective homomorphism from $F(3)$ to $F(2)$. 
Moreover it is well-known that there exist an injective homomorphism from $F(3)$ to $F(2)$. 
Suppose that $f:F(2)\to F(3)$ is a homomorphism. Let $FA(n)$ be a free abelian group of rank $n$ and $\alpha_n:F(n)\to FA(n)$ the abelianization homomorphism. Then there is a homomorphism $\hat{f}:FA(2)\to FA(3)$ such that $\alpha_3\circ f=\hat{f}\circ\alpha_2$. Then $|F(3)/f(F(2))|\geq|\alpha_3(F(3))/\alpha_3(f(F(2)))|=|\alpha_3(F(3))/(\hat{f}(\alpha_2(F(2)))|=|FA(3)/\hat{f}(FA(2))|=\infty$. Therefore $m_{\rm Coker}(f)=|F(3)/f(F(2))|=\infty$. 
After all we have $m_{\rm Ker}(F(2):F(3))=m_{\rm Ker}(F(3):F(2))=1$, $m_{\rm Coker}(F(2):F(3))=\infty$ and $m_{\rm Coker}(F(3):F(2))=1$. Therefore we have $d_{m_{\rm Ker}}(F(2),F(3))=0$ and $d_{m_{\rm Coker}}(F(2),F(3))=\infty$. 

However we have the following proposition.

\vskip 3mm

\begin{Proposition}\label{finite group}
Let $G$ and $H$ be finite groups. Then $d_{m_{\rm Ker}}(G,H)=d_{m_{\rm Coker}}(G,H)$.
\end{Proposition}

\vskip 3mm

\noindent{\bf Proof.} 
Since $d_{m_{\rm Ker}}(G,H)={\rm log}_e(m_{\rm Ker}(G:H)m_{\rm Ker}(H:G))$ and $d_{m_{\rm Coker}}(G,H)={\rm log}_e(m_{\rm Coker}(G:H)m_{\rm Coker}(H:G))$ it is sufficient to show $m_{\rm Ker}(G:H)m_{\rm Ker}(H:G)=m_{\rm Coker}(G:H)m_{\rm Coker}(H:G)$. 
Let $f:G\to H$ be a homomorphism. Since $G/{\rm Ker}(f)$ is isomorphic to $f(G)$ we have $|G|/|{\rm Ker}(f)|=|f(G)|$. On the other hand $|{\rm Coker}(f)|=|H|/|f(G)|$. Thus we have 
\[
|{\rm Ker}(f)|=\frac{|G|}{|H|}|{\rm Coker}(f)|.
\]
Therefore $|{\rm Ker}(f)|=m_{\rm Ker}(G:H)$ if and only if $|{\rm Coker}(f)|=m_{\rm Coker}(G:H)$. 
Similarly, for a homomorphism $g:G\to H$ we have
\[
|{\rm Ker}(g)|=\frac{|H|}{|G|}|{\rm Coker}(g)|
\]
and $|{\rm Ker}(g)|=m_{\rm Ker}(H:G)$ if and only if $|{\rm Coker}(g)|=m_{\rm Coker}(H:G)$. 
Thus taking such $f$ and $g$ we have
\[
m_{\rm Ker}(G:H)m_{\rm Ker}(H:G)=|{\rm Ker}(f)|\cdot|{\rm Ker}(g)|=\frac{|G|}{|H|}|{\rm Coker}(f)|\frac{|H|}{|G|}|{\rm Coker}(g)|
\]
\[
=|{\rm Coker}(f)|\cdot|{\rm Coker}(g)|=m_{\rm Coker}(G:H)m_{\rm Coker}(H:G).
\]
This completes the proof. $\Box$

\vskip 3mm

\section{Multiplicity of finitely generated modules over a principal ideal domain}\label{module} 

Throughout this section $R$ is a principal ideal domain. Some of the results in this section hold for any unitary ring $R$. See Remark \ref{unitary ring}. However we restrict our attention to a principal ideal domain for the simplicity. 
Let $M$ be an $R$-module finitely generated over $R$. We denote the minimal number of generators of $M$ over $R$ by $r(M)$. For $M=\{0\}$ we define $r(M)=0$. 
Note that when $R$ is a field, $M$ is a finite dimensional vector space over $R$ whose dimension ${\rm dim}(M)=r(M)$. 
Let $N$ be an $R$-module finitely generated over $R$ and $f:M\to N$ an $R$-linear map from $M$ to $N$. 
By ${\rm Ker}(f)$ we denote the kernel of $f$ and by ${\rm Coker}(f)=N/f(M)$ we denote the cokernel of $f$. 
We define the {\it kernel-rank-multiplicity} of $f$ by $m_{r({\rm Ker})}(f)=e^{r({\rm Ker}(f))}$ and the {\it cokernel-rank-multiplicity} of $f$ by $m_{r({\rm Coker})}(f)=e^{r({\rm Coker}(f))}$. 

\vskip 3mm

\begin{Proposition}\label{kernel and cokernel} 
(1) A kernel-rank-multiplicity $m_{r({\rm Ker})}$ is a multiplicity. 

\noindent
(2) A cokernel-rank-multiplicity $m_{r({\rm Coker})}$ is a multiplicity. 
\end{Proposition}

\vskip 3mm

For the proof of Proposition \ref{kernel and cokernel} we prepare the following proposition. 

\vskip 3mm

\begin{Proposition}\label{kernel and cokernel2} 
Let $M$, $N$ and $L$ be $R$-modules finitely generated over $R$. Let $f:M\to N$ and $g:N\to L$ be $R$-linear maps. 
Then the following (1) and (2) hold. 

\noindent
(1) $r({\rm Ker}(g\circ f))\leq r({\rm Ker}(f))+r({\rm Ker}(g))$. 

\noindent
(2) $r({\rm Coker}(g\circ f))\leq r({\rm Coker}(f))+r({\rm Coker}(g))$. 
\end{Proposition}

\vskip 3mm

For the proof of Proposition \ref{kernel and cokernel2}, we need the following lemmas. 
Since they are standard facts in module theory, we omit the proofs. See for example \cite[Chapter 3]{A-W} for the proofs.

\vskip 3mm

\begin{Lemma}\label{surjection}
Let $M$ and $N$ be $R$-modules. Suppose that $M$ is finitely generated over $R$ and there exist a surjective $R$-linear map from $M$ to $N$. Then $N$ is finitely generated over $R$ and $r(N)\leq r(M)$. 
\end{Lemma}

\vskip 3mm

\begin{Lemma}\label{submodule}
Let $M$ be an $R$-module finitely generated over $R$ and $N$ an $R$-submodule of $M$. 
Then $N$ is finitely generated over $R$ and $r(N)\leq r(M)$. 
\end{Lemma}

\vskip 3mm

\begin{Lemma}\label{kernel}
Let $M$ and $N$ be $R$-modules and $f:M\to N$ a surjective $R$-linear map. Suppose that the $R$-modules ${\rm Ker}(f)$ and $N$ are finitely generated over $R$. 
Then $M$ is finitely generated over $R$ and $r(M)\leq r({\rm Ker}(f))+r(N)$. 
\end{Lemma}

\vskip 3mm

\noindent{\bf Proof of Proposition \ref{kernel and cokernel2}.} 
(1) Since ${\rm Ker}(g\circ f)=f^{-1}({\rm Ker}(g))=f^{-1}({\rm Ker}(g)\cap f(M))$, $f|_{{\rm Ker}(g\circ f)}:{\rm Ker}(g\circ f)\to {\rm Ker}(g)\cap f(M)$ is a surjective $R$-linear map whose kernel ${\rm Ker}(f|_{{\rm Ker}(g\circ f)})$ is equal to ${\rm Ker}(f)$. By Lemma \ref{submodule} $r({\rm Ker}(g)\cap f(M))\leq r({\rm Ker}(g))$. By Lemma \ref{kernel} $r({\rm Ker}(g\circ f))\leq r({\rm Ker}(f|_{{\rm Ker}(g\circ f)}))+r({\rm Ker}(g)\cap f(M))$. Thus we have $r({\rm Ker}(g\circ f))\leq r({\rm Ker}(f))+r({\rm Ker}(g)\cap f(M))\leq r({\rm Ker}(f))+r({\rm Ker}(g))$. 

\noindent
(2) Note that ${\rm Coker}(g\circ f)=L/g(f(M))$ and $L\supset g(N)\supset g(f(M))$. By the third isomorphism theorem 
\[
(L/g(f(M)))/(g(N)/g(f(M)))\cong L/g(N).
\]
Then by Lemma \ref{kernel} $r({\rm Coker}(g\circ f))=r(L/g(f(M)))\leq r(g(N)/g(f(M)))+r(L/g(N))=r(g(N)/g(f(M)))+r({\rm Coker}(g))$. 
There exists a surjective $R$-linear map from ${\rm Coker}(f)=N/f(M)$ to $g(N)/g(f(M))$ induced by $g$. Therefore we have $r(g(N)/g(f(M)))\leq r({\rm Coker}(f))$ by Lemma \ref{surjection}. Thus we have the desired inequality. 
$\Box$ 

\vskip 3mm

\noindent{\bf Proof of Proposition \ref{kernel and cokernel}.} 
Let $M$, $N$ and $L$ be $R$-modules finitely generated over $R$ and let $f:M\to N$ and $g:N\to L$ be $R$-linear maps. 

\noindent
(1) Since $r({\rm Ker}({\rm id}_M))=r(\{0\})=0$, we have $m_{r({\rm Ker})}({\rm id}_M)=e^0=1$. 
By Proposition \ref{kernel and cokernel2} (1) we have $r({\rm Ker}(g\circ f))\leq r({\rm Ker}(f))+r({\rm Ker}(g))$. Therefore we have
\[
m_{r({\rm Ker})}(g\circ f)=e^{r({\rm Ker}(g\circ f))}\leq e^{r({\rm Ker}(f))}e^{r({\rm Ker}(g))}=m_{r({\rm Ker})}(f)m_{r({\rm Ker})}(g).
\]

\noindent
(2) Since $r({\rm Coker}({\rm id}_M))=r(\{0\})=0$, we have $m_{r({\rm Coker})}({\rm id}_M)=e^0=1$. 
By Proposition \ref{kernel and cokernel2} (2) we have $r({\rm Coker}(g\circ f))\leq r({\rm Coker}(f))+r({\rm Coker}(g))$. Therefore we have
\[
m_{r({\rm Coker})}(g\circ f)=e^{r({\rm Coker}(g\circ f))}\leq e^{r({\rm Coker}(f))}e^{r({\rm Coker}(g))}=m_{r({\rm Coker})}(f)m_{r({\rm Coker})}(g).
\]
This completes the proof. 
$\Box$ 

\vskip 3mm

By Lemma \ref{submodule} and Lemma \ref{surjection} we have both $m_{r({\rm Ker})}$ and $m_{r({\rm Coker})}$ have finite multiplicity property in the category of $R$-modules finitely generated over $R$ and $R$-linear maps. Thus we have two pseudo distances $d_{m_{r({\rm Ker})}}$ and $d_{m_{r({\rm Coker})}}$. 
The following example shows a difference between them. Let $R={\mathbb Z}$ be the ring of integers and $n{\mathbb Z}$ the ideal of ${\mathbb Z}$ generated by an integer $n$. Then both ${\mathbb Z}$ and ${\mathbb Z}/2{\mathbb Z}$ are ${\mathbb Z}$-modules finitely generated over ${\mathbb Z}$. 
It is easy to see $m_{r({\rm Ker})}({\mathbb Z}:{\mathbb Z}/2{\mathbb Z})=e$, $m_{r({\rm Ker})}({\mathbb Z}/2{\mathbb Z}:{\mathbb Z})=e$, $m_{r({\rm Coker})}({\mathbb Z}:{\mathbb Z}/2{\mathbb Z})=1$ and $m_{r({\rm Coker})}({\mathbb Z}/2{\mathbb Z}:{\mathbb Z})=e$. 
Therefore we have $d_{m_{r({\rm Ker})}}({\mathbb Z},{\mathbb Z}/2{\mathbb Z})=2$ and $d_{m_{r({\rm Coker})}}({\mathbb Z},{\mathbb Z}/2{\mathbb Z})=1$. 

In contrast to the example above, we have the following result for free $R$-modules finitely generated over $R$. 
Let $n$ be a non-negative integer. By $R^n=\underbrace{R\oplus\cdots\oplus R}_n$, we denote a free $R$-module of rank $r(R^n)=n$. Here $R^0=\{0\}$ denotes a zero-module. 

\vskip 3mm

\begin{Proposition}\label{free module} 
Let $n$ and $m$ be non-negative integers. Then 
\[
d_{m_{r({\rm Ker})}}(R^n,R^m)=d_{m_{r({\rm Coker})}}(R^n,R^m)=|m-n|. 
\]
\end{Proposition}

\vskip 3mm

\noindent{\bf Proof.} 
We may suppose without loss of generality that $m\leq n$. 
Then it is easy to see $m_{r({\rm Ker})}(R^m:R^n)=1$, $m_{r({\rm Ker})}(R^n:R^m)\leq e^{n-m}$, $m_{r({\rm Coker})}(R^m:R^n)\leq e^{n-m}$ and $m_{r({\rm Coker})}(R^n:R^m)=1$. Let $f:R^n\to R^m$ be an $R$-linear map. Then by Lemma \ref{kernel} we have $r(R^n)\leq r({\rm Ker}(f))+r(f(R^n))$. Therefore $r({\rm Ker}(f))\geq r(R^n)-r(f(R^n))$. Since $r(f(R^n))\leq r(R^m)$ by Lemma \ref{submodule} we have $r({\rm Ker}(f))\geq r(R^n)-r((R^m)=n-m$. Therefore $m_{r({\rm Ker})}(R^n:R^m)\geq e^{n-m}$. Hence $m_{r({\rm Ker})}(R^n:R^m)=e^{n-m}$. 
Let $g:R^m\to R^n$ be an $R$-linear map. Then by Lemma \ref{surjection} $r(g(R^m))\leq r(R^m)=m$. By Lemma \ref{kernel} we have $n=r(R^n)\leq r(g(R^m))+r(R^n/g(R^m))$. Therefore $r(R^n/g(R^m))\geq n-m$. Therefore $m_{r({\rm Coker})}(R^m:R^n)\geq e^{n-m}$. 
Hence $m_{r({\rm Coker})}(R^m:R^n)=e^{n-m}$. 
Then by a calculation we have the result. 
$\Box$

\vskip 3mm

The following is an immediate corollary of Proposition \ref{free module}. 

\vskip 3mm

\begin{Corollary}\label{vector space} 
Let $R$ be a field and $V$ and $W$ finite dimensional vector spaces over $R$. Then 
\[
d_{m_{r({\rm Ker})}}(V,W)=d_{m_{r({\rm Coker})}}(V,W)=|{\rm dim}(V)-{\rm dim}(W)|. 
\]
\end{Corollary}

\vskip 3mm

\begin{Theorem}\label{module1}
Let $R$ be a principal ideal domain. Let ${\mathcal O}_R$ be a set of $R$-modules finitely generated over $R$. Let ${\mathcal C}_R$ be a category whose objects are the elements of ${\mathcal O}_R$ and whose morphisms are $R$-linear maps between the elements of ${\mathcal O}_R$. 

\noindent
(1) Let $m_{r({\rm Ker})}$ be the kernel-rank-multiplicity on ${\mathcal C}_R$ and $d_{m_{r({\rm Ker})}}$ the kernel-rank-multiplicity distance. Let $M$ and $N$ be elements of ${\mathcal O}_R$. Then $M$ and $N$ are isomorphic if and only if $d_{m_{r({\rm Ker})}}(M,N)=0$. Thus the pseudo distance $d_{m_{r({\rm Ker})}}$ defines a distance on the set of isomorphism classes of ${\mathcal O}_R$. 

\noindent
(2) Let $m_{r({\rm Coker})}$ be the cokernel-rank-multiplicity on ${\mathcal C}_R$ and $d_{m_{r({\rm Coker})}}$ the cokernel-rank-multiplicity distance. Let $M$ and $N$ be elements of ${\mathcal O}_R$. Then $M$ and $N$ are isomorphic if and only if $d_{m_{r({\rm Coker})}}(M,N)=0$. Thus the pseudo distance $d_{m_{r({\rm Coker})}}$ defines a distance on the set of isomorphism classes of ${\mathcal O}_R$. 
\end{Theorem}

\vskip 3mm

The essential part of Theorem \ref{module1} is the following proposition. 

\vskip 3mm

\begin{Proposition}\label{module2}
Let $M$ and $N$ be $R$-modules finitely generated over $R$. 

\noindent
(1) Suppose that there exist injective $R$-linear maps $f:M\to N$ and $g:N\to M$. Then $M$ and $N$ are isomorphic. 

\noindent
(2) Suppose that there exist surjective $R$-linear maps $f:M\to N$ and $g:N\to M$. Then both $f$ and $g$ are injective, and $M$ and $N$ are isomorphic. 

\end{Proposition}

\vskip 3mm

For the proof of Proposition \ref{module2}, we prepare some lemmas. For an $R$-module $M$, we denote the torsion submodule of $M$ by ${\rm tor}(M)$. 
\vskip 3mm

\begin{Lemma}\label{torsion}
Let $M$ and $N$ be $R$-modules and $f:M\to N$ an $R$-linear map. Then $f({\rm tor}(M))\subset{\rm tor}(N)$. 
\end{Lemma}

\vskip 3mm

\noindent{\bf Proof.} 
Suppose $x\in{\rm tor}(M)$. Then there is a non-zero element $r\in R$ such that $rx=0$. Then $rf(x)=f(rx)=f(0)=0$. Therefore $f(x)\in{\rm tor}(N)$. 
$\Box$

\vskip 3mm

Thus for $f:M\to N$, an induced $R$-linear map $\hat{f}:M/{\rm tor}(M)\to N/{\rm tor}(N)$ is defined by $\hat{f}([x])=[f(x)]$ where $x\in M$ and $[x]\in M/{\rm tor}(M)$ denotes the set of elements $y\in M$ with $x-y\in{\rm tor}(M)$. 

\vskip 3mm

\begin{Lemma}\label{torsion2}
Let $M$ and $N$ be $R$-modules and $f:M\to N$ an injective $R$-linear map. Let $\hat{f}:M/{\rm tor}(M)\to N/{\rm tor}(N)$ be an induced $R$-linear map. 
Then $\hat{f}$ is injective. 
\end{Lemma}

\vskip 3mm

\noindent{\bf Proof.} 
Suppose $\hat{f}([x])=\hat{f}([y])$. Then $[f(x)]=[f(y)]$ and therefore $f(x)-f(y)=f(x-y)\in{\rm tor}(N)$. Let $r\in R$ such that $rf(x-y)=0$. Since $rf(x-y)=f(r(x-y))$ and $f$ is injective we have $r(x-y)=0$. Thus $x-y\in{\rm tor}(M)$ and $[x]=[y]$. 
$\Box$

\vskip 3mm

For an element $r$ of $R$, we denote the ideal generated by $r$ by $(r)$. Then the quotient $R/(r)$ is a torsion $R$-module unless $r=0$. 
The following Lemma \ref{Artinian} and Proposition \ref{Artinian2} are known facts. However we give a direct elementary proof of Lemma \ref{Artinian} for convenience.

\vskip 3mm

\begin{Lemma}\label{Artinian}
Let $p$ be a prime element of $R$ and $k$ a natural number. Then an $R$-module $R/(p^k)$ is Artinian. 
\end{Lemma}

\vskip 3mm

\noindent{\bf Proof.} 
It is sufficient to show that a proper $R$-submodule $M$ of $R/(p^k)$ is isomorphic to $R/(p^j)$ for some non-negative integer $j<k$. 
Let $\varphi:R\to R/(p^k)$ be a canonical projection. Then $\varphi^{-1}(M)$ is an $R$-submodule, hence an ideal of $R$. Let $a$ be an element of $R$ such that $\varphi^{-1}(M)=(a)$. Since $(a)\supset(p^k)$ we have $(a)=(p^l)$ for some non-negative integer $l\leq k$. If $l=0$ then $\varphi^{-1}(M)=(1_R)=R$ where $1_R$ is the multiplicative identity of $R$ and then $M=\varphi(R)=R/(p^k)$ is not a proper $R$-submodule. Thus we have $l\geq1$. Let $\psi:R\to R$ be an $R$-linear map defined by $\psi(x)=p^lx$. Then $\psi$ induces an $R$-linear map $\hat{\psi}:R/(p^{k-l})\to R/(p^k)$. By the fact that $R$ is a unique factorization domain it follows that $\psi$ is injective. By the definition of $\hat{\psi}$ we have $\hat{\psi}(R/(p^{k-l}))=M$. Therefore $M$ is isomorphic to $R/(p^{k-l})$. 
$\Box$

\vskip 3mm

It is well-known that a direct sum of finitely many Artinian $R$-modules is Artinian. By the structure theorem, a torsion $R$-module finitely generated over $R$ is a direct sum of finitely many cyclic $R$-modules of the form $R/(p^k)$. Therefore we have the following proposition. 

\vskip 3mm

\begin{Proposition}\label{Artinian2}
A torsion $R$-module finitely generated over $R$ is Artinian. 
\end{Proposition}

\vskip 3mm

\noindent{\bf Proof of Proposition \ref{module2}.} 
(1) By Lemma \ref{torsion2} the free $R$-modules $M/{\rm tor}(M)$ and $N/{\rm tor}(N)$ have the same rank. Therefore they are isomorphic. By Lemma \ref{torsion} we have injective $R$-linear maps $f':{\rm tor}(M)\to{\rm tor}(N)$ and $g':{\rm tor}(N)\to{\rm tor}(M)$. Then $g'\circ f':{\rm tor}(M)\to{\rm tor}(M)$ is an injective $R$-linear map. Since ${\rm tor}(M)$ is Artinian, it is co-Hopfian. Therefore $g'\circ f'$ is surjective. Thus $g'$ is a surjection, hence an isomorphism. Then by the structure theorem of $R$-modules finitely generated over $R$, $M$ and $N$ are isomorphic. 

\noindent
(2) Note that $g\circ f:M\to M$ is a surjective $R$-linear map. It is well-known that $M$ is Noetherian, hence Hopfian. Therefore $g\circ f$ is injective. Therefore $f$ is an injection, hence an isomorphism. Similarly, $g$ is also an isomorphism. 
$\Box$

\vskip 3mm

\noindent{\bf Proof of Theorem \ref{module1}.} 
(1) The \lq only if' part is clear. Suppose $d_{m_{r({\rm Ker})}}(M,N)=0$. Then there are injective $R$-linear maps $f:M\to N$ and $g:N\to M$. Then by Proposition \ref{module2} (1) $M$ and $N$ are isomorphic. 

\noindent
(2) The \lq only if' part is clear. Suppose $d_{m_{r({\rm Coker})}}(M,N)=0$. Then there are surjective $R$-linear maps $f:M\to N$ and $g:N\to M$. Then by Proposition \ref{module2} (2) $M$ and $N$ are isomorphic. 
$\Box$

\vskip 3mm

\begin{Remark}\label{unitary ring}
{\rm The main reason that we have restricted our attention to the case $R$ is a principal ideal domain is the fact that Lemma \ref{submodule} does not hold in general when $R$ is not a principal ideal domain. Let $R={\mathbb Z}[X]$ be a polynomial ring over the integer ring ${\mathbb Z}$. 
Let $M=R$ and $N=(2,X)$ where $(2,X)$ denotes the ideal of $R$ generated by $2$ and $X$. Then $r(M)=1$. 
Since $N$ is an ideal of an $R$-module $M$, $N$ is an $R$-submodule of $M$. Since $N$ is not a principal ideal we have $r(N)=2$. 
Actually Proposition \ref{kernel and cokernel2} (1) does not hold in general when $R$ is not a principal ideal domain. 
Let $L=\{0\}$ be the zero module. Let $f:N\to M$ be the inclusion map and $g:M\to L$ the zero map. 
Then $r({\rm Ker}(g\circ f))=r(N)=2$, $r({\rm Ker}(f))=r(\{0\})=0$ and $r({\rm Ker}(g))=r(M)=1$. 
Moreover Proposition \ref{module2} (1) does not hold in general when $R$ is not a principal ideal domain. 
Let $f:M\to N$ be an $R$-linear map defined by $f(x)=2x$ and $g:N\to M$ the inclusion $R$-linear map. Then they are injective. 
Since $r(M)\neq r(N)$ they are not isomorphic. 

We note however Lemma \ref{surjection} and Lemma \ref{kernel} hold for any unitary ring $R$. Therefore Proposition \ref{kernel and cokernel2} (2) and Proposition \ref{kernel and cokernel} (2) also hold for any unitary ring $R$. 

We also note that Proposition \ref{module2} (2) is true when $M$ or $N$ is a Hopfian $R$-module. 
}
\end{Remark}

\vskip 3mm

\section{Multiplicity of knots}\label{knot} 

\vskip 3mm

In this section we define a category of oriented knot types in $3$-space and define a multiplicity on it. Then we study its multiplicity distance. 
Throughout this section we work in the piecewise linear category. 
First we define a category which we call a {\it neighbourhood category of oriented knots}. 

For a manifold $X$ with boundary, we denote the interior of $X$ by ${\rm int}(X)$ and the boundary of $X$ by $\partial(X)$ respectively. 
By ${\mathbb D}^2$ we denote the unit disk in the $2$-dimensional Euclidean space ${\mathbb R}^2$ and by ${\mathbb S}^1=\partial({\mathbb D}^2)$ we denote the unit circle in ${\mathbb R}^2$. 
A {\it core} of a solid torus $V$ is a circle $c\subset{\rm int}(V)$ such that the pair $(V,c)$ is homeomorphic to the pair $({\mathbb D}^2\times{\mathbb S}^1,(0,0)\times{\mathbb S}^1)$. 
Two oriented cores $c_1$ and $c_2$ of $V$ are {\it equivalent} if they are ambient isotopic in $V$ as oriented circles. A {\it core-orientation} of $V$ is an equivalence class of oriented cores in $V$. A solid torus $V$ is {\it core-oriented} if a core-orientation of $V$ is specified. 
Let $V$ be a core-oriented solid torus embedded in $\mathbb{S}^3$ and $k$ an oriented circle contained in ${\rm int}(V)$. 
Then the pair $(V,k)$ is called a {\it nesting}. Two nestings $(V,k)$ and $(V',k')$ are {\it equivalent} if there exist an orientation preserving homeomorphism $\varphi:{\mathbb S}^3\to{\mathbb S}^3$ such that $\varphi(V)=V'$ respecting core-orientations and $\varphi(k)=k'$ respecting orientations. 
The equivalence class containing $(V,k)$ is denoted by $[(V,k)]$. 
Let ${\mathcal K}$ be the set of all oriented knot types in the $3$-sphere ${\mathbb S}^3$. 
Let $K_1$ and $K_2$ be elements of ${\mathcal K}$. Namely $K_1$ and $K_2$ are oriented knot types in $\mathbb{S}^3$. 
A {\it morphism} $f$ from $K_1$ to $K_2$ is an equivalence class $[(V,k)]$ of nestings such that the knot type of an oriented core of $V$ is $K_2$ and the knot type of $k$ is $K_1$. Then we denote $f:K_1\to K_2$ and $f=[(V,k)]$. 

Let $K_3$ be an element of ${\mathcal K}$ and $g=[(W,l)]$ a morphism from $K_2$ to $K_3$. 
Let $N\subset{\rm int}(W)$ be a regular neighbourhood of $l$. Then $N$ is a solid torus and $l$ is a core of $N$. We give a core-orientation of $N$ by the orientation of $l$. Let $f:{\mathbb S}^3\to{\mathbb S}^3$ be an orientation preserving homeomorphism such that $f(V)=N$ respecting core-orientations. Then $f(k)$ is an oriented knot contained in $W$ whose knot type is $K_1$. Thus $[(W,f(k))]$ is a morphism from $K_1$ to $K_3$. We define the composition of $f$ and $g$ by $g\circ f=[(W,f(k))]$. 

We will check the well-definedness of composition as follows. 
Suppose that an orientation preserving homeomorphism $\varphi:{\mathbb S}^3\to{\mathbb S}^3$ maps $(V,k)$ to $(V',k')$ and an orientation preserving homeomorphism $\psi:{\mathbb S}^3\to{\mathbb S}^3$ maps $(W,l)$ to $(W',l')$. 
Let $N'\subset{\rm int}(W')$ be a regular neighbourhood of $l'$ and let $f':{\mathbb S}^3\to{\mathbb S}^3$ be an orientation preserving homeomorphism such that $f(V')=N'$ respecting core-orientations. 
It is sufficient to show $[(W,f(k))]=[(W',f'(k'))]$. Since $[(W,f(k))]=[W',\psi\circ f(k))]$ by $\psi$, it is sufficient to show that $\psi\circ f(k)$ is ambient isotopic to $f'(k')$ in $W'$. 
Note that both $\psi(N)$ and $N'$ are solid tori in ${\rm int}(W')$ with the same oriented core $l'$. Therefore there is an isotopy $h_t:W'\to W'$ with $0\leq t\leq1$ and $h_0={\rm id}_{W'}$ point-wisely fixed on $\partial(W')\cup l'$ such that $h_1(\psi(N))=N'$. Let $\gamma:{\mathbb S}^3\to{\mathbb S}^3$ be an extension of $h_1:W'\to W'$ such that $\gamma(x)=x$ for any $x\in{\mathbb S}^3\setminus W'$. 
Then $\psi\circ f(k)$ is ambient isotopic to $h_1(\psi\circ f(k))=\gamma(\psi\circ f(k))=\gamma\circ\psi\circ f(k)$. 
Thus it is sufficient to show that $\gamma\circ\psi\circ f(k)$ is ambient isotopic to $f'(k')$ in $W'$. 
We now have the following sequence of orientation preserving self-homeomorphisms of ${\mathbb S}^3$. 
\[
(N',f'(k'))\stackrel{f'}{\leftarrow}(V',k')\stackrel{\varphi}{\leftarrow}(V,k)\stackrel{f}{\rightarrow}(N,f(k))\stackrel{\psi}{\rightarrow}(\psi(N),\psi\circ f(k))\stackrel{\gamma}{\rightarrow}(N',\gamma\circ\psi\circ f(k))
\]
Let $\eta:N'\to N'$ be the restriction of $\gamma\circ\psi\circ f\circ\varphi^{-1}\circ(f')^{-1}$ to $N'$. Then $\eta$ is an orientation preserving homeomorphism also preserving the core-orientation of $N'$. Since $\eta$ is a restriction of a self-homeomorphism of ${\mathbb S}^3$, $\eta$ maps a preferred longitude of $N'$ to a preferred longitude of $N'$. Therefore $\eta$ is isotopic to ${\rm id}_{N'}$, see for example \cite{Rolfsen}. Therefore $f'(k')$ is ambient isotopic to $\eta(f'(k'))=\gamma\circ\psi\circ f(k)$ in $N'$, hence in $W'$ as desired. 

\begin{figure}[htbp]
      \begin{center}
\scalebox{0.5}{\includegraphics*{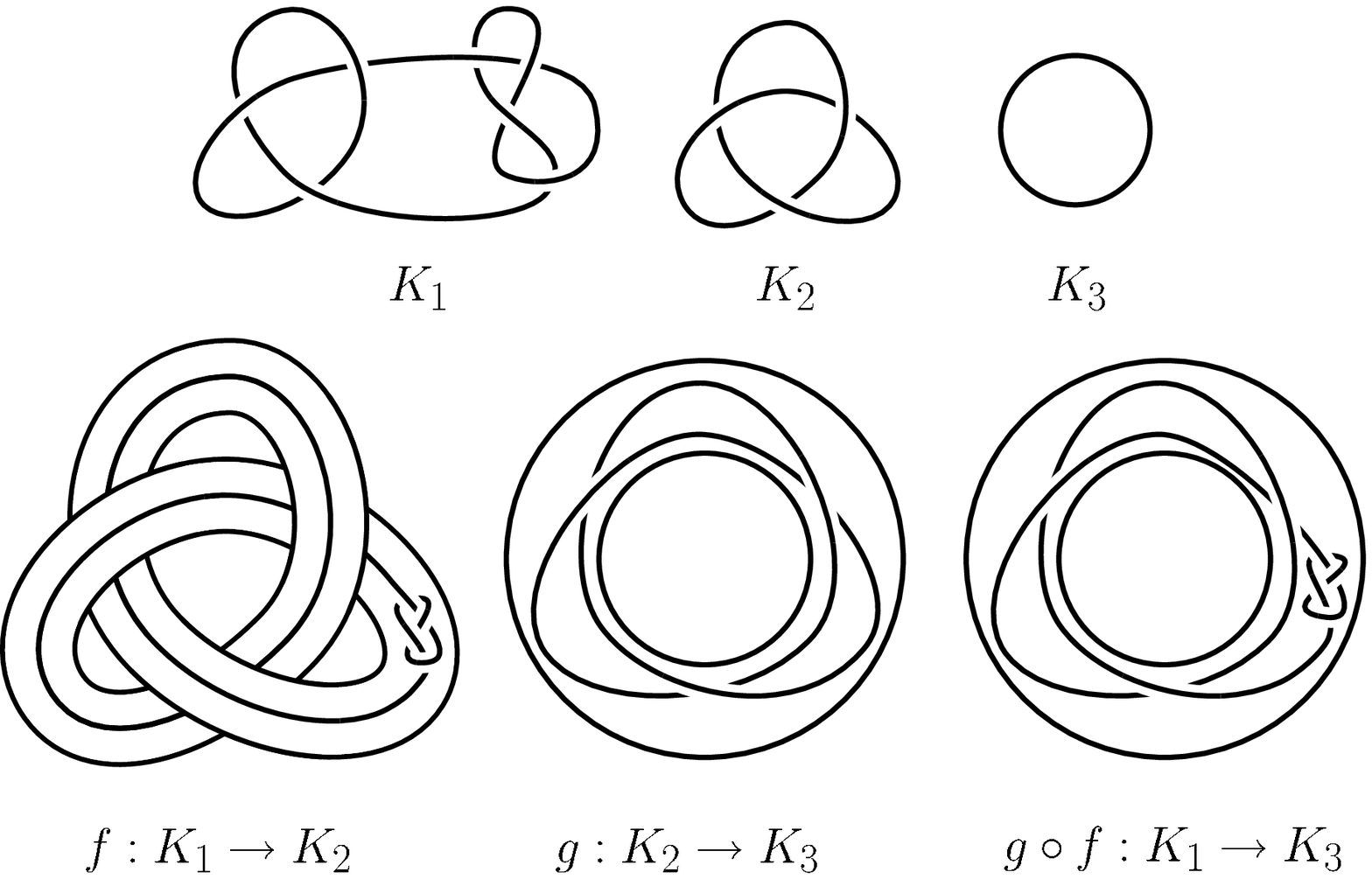}}
      \end{center}
   \caption{}
  \label{knot-morphism}
\end{figure}

An example of the composition of two morphisms is illustrated in Figure \ref{knot-morphism}. 

Let $K$ be an element of ${\mathcal K}$. Let $k\subset{\mathbb S}^3$ be an oriented knot whose knot type is $K$ and $N$ a regular neighbourhood of $k$ in ${\mathbb S}^3$. We define ${\rm id}_K:K\to K$ by ${\rm id}_K=[(N,k)]$. It is clear that ${\rm id}_K$ actually satisfies the axiom of identity morphism of category. By the definition other axioms of category clearly hold. Thus we have a category whose objects are the elements of ${\mathcal K}$. We call this category a {\it neighbourhood category of oriented knots} and denote it by ${\mathcal C}_{\mathcal K}$.

\vskip 3mm

\begin{Remark}\label{exterior}
{\rm In \cite{Nikkuni} Ryo Nikkuni has pointed out the following. Suppose that $f=[(V,k)]$ is a morphism from $K_1$ to $K_2$. Let $N\subset{\rm int}(V)$ be a regular neighbourhood of $k$. Let $E_2={\mathbb S}^3\setminus{\rm int}(V)$ and $E_1={\mathbb S}^3\setminus{\rm int}(N)$. Then $E_i$ is an exterior manifold of a knot representing $K_i$ for $i=1,2$ and we have an inclusion $E_2\subset E_1$. Then the inclusion induces a homomorphism from $\pi_1(E_2)$ to $\pi_1(E_1)$ where $\pi_1(X)$ denotes the fundamental group of $X$. 
}
\end{Remark}

\vskip 3mm

For any elements $K_1$ and $K_2$ of ${\mathcal K}$, there is a morphism $c:K_1\to K_2$ with $c=[(N,k)]$ where $N$ is a core-oriented solid torus in ${\mathbb S}^3$ whose knot type is $K_2$ and $k$ is an oriented knot contained in a $3$-ball in $N$ whose knot type is $K_1$. 
We call $c$ the {\it null morphism} from $K_1$ to $K_2$. By definition it is clear that both $c\circ f$ and $g\circ c$ are again null morphisms for any $f$ and $g$. For the existence of non-null morphisms, we have the following proposition. For an oriented knot $k$ in ${\mathbb S}^3$, $-k$ denotes an oriented knot obtained from $k$ by reversing its orientation. For an element $K$ of ${\mathcal K}$ represented by an oriented knot $k$ in ${\mathbb S}^3$, $-K$ denotes an element of ${\mathcal K}$ represented by $-k$. A morphism $f=[(V,k)]$ from $K_1$ to $K_2$ is called an {\it inversion} if $k$ is a core of $V$ and the core-orientation of $V$ is represented by $-k$. Thus if $f:K_1\to K_2$ is an inversion then $K_2=-K_1$. 

\vskip 3mm

\begin{Proposition}\label{satellite}
Let $K_1$ and $K_2$ be elements of ${\mathcal K}$. Suppose that there exists a non-null morphism $f:K_1\to K_2$. 
Then at least one of the following holds. 

\noindent
(1) $K_1=K_2$ and $f={\rm id}_{K_1}$. 

\noindent
(2) $K_1=-K_2$ and $f$ is an inversion.

\noindent
(3) $K_2$ is the trivial knot type. 

\noindent
(4) $K_1$ is a satellite knot type with companion knot type $K_2$. 
\end{Proposition}

\vskip 3mm

\noindent{\bf Proof.} 
Suppose $f=[(V,k)]$. Suppose that $K_2$ is non-trivial. Then $V$ is a knotted solid torus in ${\mathbb S}^3$. Since $k$ is not contained in any $3$-ball in $V$, $k$ is essential in $V$. Then either $k$ is a core of $V$, then we have (1) or (2), or $k$ is a satellite knot with companion solid torus $V$, then we have (4). 
$\Box$ 

\vskip 3mm

Let $\varphi:X\to Y$ be a continuous map from a circle $X$ to a $1$-dimensional manifold $Y$ without boundary. 
A point $x\in X$ is a {\it regular point} of $\varphi$ if there is a neighbourhood $U$ of $x$ in $X$ such that $\varphi|_U$ is injective. A point $x\in X$ is a {\it critical point} of $\varphi$ if it is not a regular point. By ${\rm cr}(\varphi)$ we denote the set of all critical points of $\varphi$. 
We say that $\varphi$ is {\it simple} if $X=I_1\cup\cdots\cup I_n$ for some $n\in{\mathbb N}$ and closed intervals $I_1,\cdots,I_n$ such that $\varphi|_{I_i}$ is injective for each $i\in\{1,\cdots,n\}$. We say that $\varphi$ is {\it generic} if it is simple and $\varphi|_{{\rm cr}(\varphi)}$ is injective. 

We now define a multiplicity on the neighbourhood category of oriented knots as follows. 
Let $K_1$ and $K_2$ be elements of ${\mathcal K}$ and  and $f=[(V,k)]$ a morphism from $K_1$ to $K_2$. 
Let $h:V\to\mathbb{S}^1\times\mathbb{D}^2$ be a homeomorphism and $\pi:\mathbb{S}^1\times\mathbb{D}^2\to \mathbb{S}^1$ a natural projection. 
We say that $h$ is {\it generic} with respect to $k$ if $\pi\circ h|_k:k\to\mathbb{S}^1$ is generic. 
Let ${m(h,k)}={\rm max}\{|(\pi\circ h|_{k})^{-1}(y)|\ |\ y\in\mathbb{S}^1\}$. 
Note that if $y\in\pi\circ h({\rm cr}(\pi\circ h|_{k}))$ then $|(\pi\circ h|_{k})^{-1}(y)|<m(h,k)$. Then we define the multiplicity of $f$, denoted by $m_{\mathcal K}(f)$, to be the minimum of $m(h,k)$ where $h:V\to\mathbb{S}^1\times\mathbb{D}^2$ varies over all homeomorphisms that are generic with respect to $k$. By definition $m_{\mathcal K}(f)$ is well-defined. We call $m_{\mathcal K}(f)$ the {\it knot-multiplicity} of $f$. 

\vskip 3mm

\begin{Proposition}\label{knot-multiplicity}
The knot-multiplicity $m_{\mathcal K}$ is a multiplicity on the neighbourhood category of oriented knots. 
\end{Proposition}

\vskip 3mm

\begin{figure}[htbp]
      \begin{center}
\scalebox{0.5}{\includegraphics*{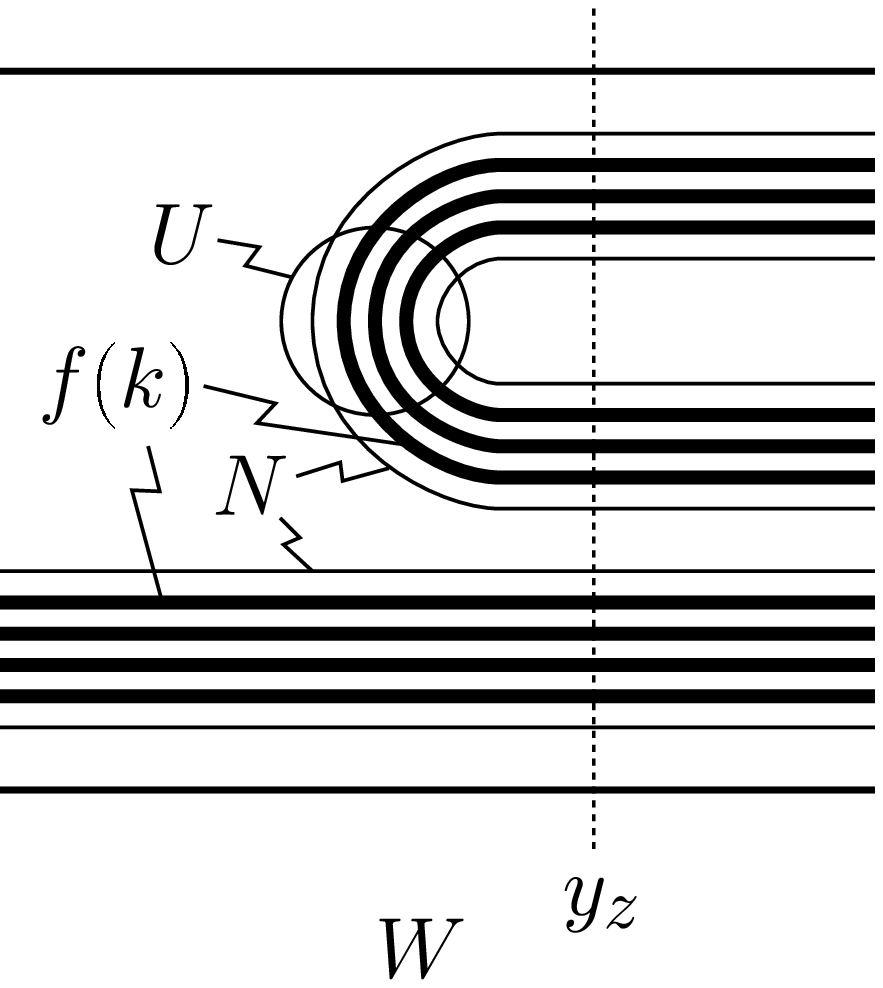}}
      \end{center}
   \caption{}
  \label{critical-point}
\end{figure}

\noindent{\bf Proof.} 
It is clear that $m_{\mathcal K}({\rm id}_K)=1$ for any element $K$ of ${\mathcal K}$. 
Let $f:K_1\to K_2$ and $g:K_2\to K_3$ be morphisms where $f=[(V,k)]$ and $g=[(W,l)]$. 
Let $h:V\to\mathbb{S}^1\times\mathbb{D}^2$ be a homeomorphism that is generic with respect to $k$ with $m(h,k)=m_{\mathcal K}(f)$ and $s:W\to\mathbb{S}^1\times\mathbb{D}^2$ a homeomorphism that is generic with respect to $l$ with $m(s,l)=m_{\mathcal K}(g)$. 
Let $N\subset{\rm int}(W)$ be a regular neighbourhood of $l$ and $f:{\mathbb S}^3\to{\mathbb S}^3$ an orientation preserving homeomorphism such that $f(V)=N$ respecting core-orientations. Then the composition $g\circ f:K_1\to K_3$ is given by $g\circ f=[(W,f(k))]$. By a modification of $s$ on $N$ we may assume that $f({\rm cr}(\pi\circ h|_k))$ is a subset of ${\rm cr}(\pi\circ s|_{f(k)})$ and away from a small neighbourhood, say $U$, of ${\rm cr}(\pi\circ s|_l)$, and all other critical points ${\rm cr}(\pi\circ s|_{f(k)})\setminus f({\rm cr}(\pi\circ h|_k))$ are contained in $U$. See Figure \ref{critical-point} which illustrates a part of $W$ containing a critical point, say $z$, of $\pi\circ s|_l$. Let $y_z$ be a point in ${\mathbb S}^1$ as indicated in Figure \ref{critical-point}. Then $|(\pi\circ s|_{f(k)})^{-1}(y_z)|$ is locally maximal in a small neighbourhood of $y_z$ in ${\mathbb S}^1$. 
Therefore we only need to consider $|(\pi\circ s|_{f(k)})^{-1}(y_z)|$ for each critical point $z$ of $\pi\circ s|_l$. 
Therefore we have $m(s,f(k))\leq m(h,k)m(s,l)$. Thus $m_{\mathcal K}(g\circ f)\leq m(s,f(k))\leq m_{\mathcal K}(f)m_{\mathcal K}(g)$ as desired. 
$\Box$

\begin{Proposition}\label{one or two}
Let $K_1$ and $K_2$ be elements of ${\mathcal K}$. 

\noindent
(1) $m_{\mathcal K}(K_1:K_2)=1$ if and only if $K_1=K_2$ or $K_1=-K_2$. 

\noindent
(2) $m_{\mathcal K}(K_1:K_2)=2$ if and only if $K_1$ is trivial and $K_2$ is non-trivial, or $K_1$ is non-trivial and $K_1$ is a $(2,p)$-cable knot of $K_2$ or $-K_2$ for some odd number $p$. (When $K_2$ is trivial, $K_1$ is a non-trivial $(2,p)$-torus knot.)
\end{Proposition}

\vskip 3mm

\noindent{\bf Proof.} 
(1) Suppose $m_{\mathcal K}(K_1:K_2)=1$. Then there exists a morphism $f=[(V,k)]$ from $K_1$ to $K_2$ with $m_{\mathcal K}(f)=1$. Then $k$ is a core of $V$ and hence we have $K_1=K_2$ or $K_1=-K_2$. The converse is similar. 

\noindent
(2) Suppose $m_{\mathcal K}(K_1:K_2)=2$. Then there exists a morphism $f=[(V,k)]$ from $K_1$ to $K_2$ with $m_{\mathcal K}(f)=2$. Then the number of critical points of $k$ with respect to the projection from $V$ to ${\mathbb S}^1$ is $0$ or $2$. In the first case $k$ is a $(2,p)$-cabling of the core of $V$. Therefore $K_1$ is a $(2,p)$-cable knot of $K_2$ or $-K_2$ for some odd number $p$. If $K_2$ is non-trivial, then $K_1$ is also non-trivial. If $K_2$ is trivial, $m_{\mathcal K}(K_1:K_2)\neq1$ implies $K_1$ is non-trivial. In the second case $k$ is a trivial knot. Then $m_{\mathcal K}(K_1:K_2)\neq1$ implies $K_2$ is non-trivial. 
The converse is clear. 
$\Box$

\vskip 3mm

\begin{Theorem}\label{unoriented knot}
Let ${\mathcal C}_{\mathcal K}$ be a neighbourhood category of oriented knots. 
Let $m_{\mathcal K}$ be the knot-multiplicity on ${\mathcal C}_{\mathcal K}$ and $d_{m_{\mathcal K}}$ the knot-multiplicity distance. Let $K_1$ and $K_2$ be elements of ${\mathcal K}$. Then $d_{m_{\mathcal K}}(K_1,K_2)=0$ if and only if $K_1=K_2$ or $K_1=-K_2$. 
Thus the pseudo distance $d_{m_{\mathcal K}}$ defines a distance on the set of unoriented knot types in ${\mathbb S}^3$. 
\end{Theorem}

\vskip 3mm

\noindent{\bf Proof.} 
Note that $d_{m_{\mathcal K}}(K_1,K_2)=0$ if and only if $m_{\mathcal K}(K_1:K_2)=m_{\mathcal K}(K_2:K_1)=1$. 
Then by Proposition \ref{one or two} (1) we have the conclusion. 
$\Box$

\vskip 3mm

We now discuss several relations between knot-multiplicity and some geometric knot invariants as follows. 
In \cite{Ozawa} Ozawa defined a {\it trunk} of a knot as follows. 
Let $k$ be an oriented knot in $\mathbb{S}^3\setminus\{(0,0,0,1),(0,0,0,-1)\}$. 
Let $\varphi:\mathbb{S}^3\setminus\{(0,0,0,1),(0,0,0,-1)\}\to\mathbb{S}^2\times\mathbb{R}$ be a fixed homeomorphism and 
$\pi:\mathbb{S}^2\times\mathbb{R}\to\mathbb{R}$ a natural projection. 
We say that $k$ is in {\it general position} if $\pi\circ\varphi|_k:k\to{\mathbb R}$ is generic. 
Then the trunk of $k$, denoted by ${\rm trunk}(k)$, is defined by ${{\rm trunk}(k)}={\rm max}\{|(\pi\circ\varphi|_k)^{-1}(y)||y\in\mathbb{R}\}$. 
The trunk of an oriented knot type $K$, denoted by ${\rm trunk}(K)$, is defined to be the minimum of ${\rm trunk}(k)$ where $k$ varies over all oriented knots in ${\mathbb S}^3$ in general position whose oriented knot type is $K$. 

\vskip 3mm

\begin{Proposition}\label{trunk}
Let $K_1$ and $K_2$ be elements of ${\mathcal K}$. Then 
\[
\frac{{\rm trunk}(K_1)}{{\rm trunk}(K_2)}\leq m_{\mathcal K}(K_1:K_2)\leq{\rm trunk}(K_1).
\]
If $K_2$ is non-trivial, $K_2\neq K_1$, $K_2\neq-K_1$ and not a companion of $K_1$, then $m_{\mathcal K}(K_1:K_2)={\rm trunk}(K_1)$. 
\end{Proposition}

\vskip 3mm

\noindent{\bf Proof.} 
Let $k_2$ be an oriented knot in ${\mathbb S}^3$ representing $K_2$ with ${\rm trunk}(k_2)={\rm trunk}(K_2)$. 
Let $N$ be a regular neighbourhood of $k_2$ and $f=[(N,k_1)]$ a morphism from $K_1$ to $K_2$ with $m_{\mathcal K}(f)=m_{\mathcal K}(K_1:K_2)$. Then by deforming $k_1$ if necessary, we have ${\rm trunk}(k_1)\leq{\rm trunk}(k_2)m_{\mathcal K}(f)$. Thus we have ${\rm trunk}(K_1)\leq{\rm trunk}(K_2)m_{\mathcal K}(K_1:K_2)$. Therefore the first inequality holds. 
Let $c=[(V,k)]$ be a null morphism from $K_1$ to $K_2$. Clearly $m_{\mathcal K}(c)\leq{\rm trunk}(K_1)$. Therefore the second inequality holds. Conversely, since $k$ is contained in a $3$-ball in $V$, $k$ lifts to the universal covering space $\tilde{V}$ of $V$. Then the projection from $V$ to ${\mathbb S}^1$ lifts to the projection from $\tilde{V}$ to ${\mathbb R}$. Therefore we have ${\rm trunk}(K_1)\leq m_{\mathcal K}(c)$. 
Then the final statement follows by Proposition \ref{satellite}. 
$\Box$

\vskip 3mm

As we have observed, there are not so many morphisms when the target knot is non-trivial. Therefore,  for an element $K$ of ${\mathcal K}$, we define $m(K)=m_{\mathcal K}(K:T)$ where $T$ is the trivial knot type, and call it the {\it multiplicity index} of $K$. Since ${\rm trunk}(T)=2$, we have the following corollary of Proposition \ref{trunk}. 

\vskip 3mm

\begin{Corollary}\label{trunk2}
Let $K$ be an element of ${\mathcal K}$. Then $\displaystyle{m(K)\geq\frac{{\rm trunk}(K)}{2}}$.
\end{Corollary}

\vskip 3mm

If $K\neq T$ then $m_{\mathcal K}(T:K)={\rm trunk}(T)=2$. Therefore 
\[
d_{m_{\mathcal K}}(K,T)\geq{\rm log}_e\frac{{\rm trunk}(K)}{2}\cdot2={\rm log}_e{\rm trunk}(K).
\]
In \cite{Ozawa} it is shown that for any natural number $n$, there exist an oriented knot type $K$ with ${\rm trunk}(K)\geq n$. 
Therefore we have the following proposition. 

\vskip 3mm

\begin{Proposition}\label{unbounded}
The pseudo metric space $(\mathcal{K},d_{m_{\mathcal K}})$ is unbounded. 
\end{Proposition}

\vskip 3mm

Let ${\rm braid}(K)$ be the braid index of $K$ and ${\rm bridge}(K)$ the bridge index of $K$. 

\vskip 3mm

\begin{Proposition}\label{braid and bridge}
Let $K$ be an element of ${\mathcal K}$. 

\noindent
(1) $m(K)\leq{\rm braid}(K)$. 

\noindent
(2) $m(K)\leq2\cdot{\rm bridge}(K)-1$. 
\end{Proposition}

\vskip 3mm

\noindent{\bf Proof.} 
(1) Suppose ${\rm braid}(K)=n$. Then $K$ has a representative $k$ contained in an unknotted solid torus $V\subset{\mathbb S}^3$ such that $k$ intersects each meridian disk of $V$ at exactly $n$ points. Then $f=[(V,k)]$ is a morphism from $K$ to $T$ with $m_{\mathcal K}(f)=n$. 

\noindent
(2) Let $B_1$ and $B_2$ be $3$-balls such that ${\mathbb S}^3=B_1\cup B_2$ and $B_1\cap B_2=\partial(B_1)=\partial(B_2)$. 
Suppose ${\rm bridge}(K)=n$. Then there is a representative $k$ of $K$ such that both $(B_1,B_1\cap k)$ and $(B_2,B_2\cap k)$ are trivial $n$-string tangles. Choose a point $x\in\partial(B_1)\cap k$. Let $\gamma$ be a circle on $\partial(B_1)\setminus k$ separating $x$ from other $2n-1$ points of $\partial(B_1)\cap k$. Let $N$ be a regular neighbourhood of $\gamma$ in ${\mathbb S}^3$ with $N\cap k=\emptyset$ intersecting $\partial(B_1)$ in an annulus. Let $V={\mathbb S}^3\setminus{\rm int}(N)$. Then $V$ is an unknotted solid torus containing $k$. Note that $V\cap\partial(B_1)$ is a disjoint union of two meridian disks of $V$, one is intersecting $k$ at $x$, and the other is intersecting $k$ at $2n-1$ points. Then there is a homeomorphism $h:V\to{\mathbb S}^1\times{\mathbb D}^2$ such that $\pi\circ h:V\to{\mathbb S}^1$ maps each of these meridian disks to a point and each of $k\cap B_1$ and $k\cap B_2$ have exactly $n-1$ critical points of $\pi\circ f|_k$. Therefore $f=[(V,k)]$ is a morphism from $K$ to $T$ with $m_{\mathcal K}(f)\leq2n-1$. 
$\Box$

\vskip 3mm

We now consider knot types $K$ with $m(K)\leq n$ for $n=1,2,3,4$ as follows. Note that the condition $m(K)\leq n$ is equivalent to the condition $d_{m_{\mathcal K}}(K,T)\leq{\rm log}_e2n$.

\vskip 3mm

\begin{Proposition}\label{multiplicity index}
Let $K$ be an element of ${\mathcal K}$. 

\noindent
(1) $m(K)=1$ if and only if $K$ is trivial. 

\noindent
(2) $m(K)=2$ if and only if $K$ is a $(2,p)$-torus knot for some odd number $p\neq\pm1$. 

\noindent
(3) $m(K)=3$ if and only if ${\rm braid}(K)=3$ or $K$ is a connected sum of some $2$-bridge knots. 

\noindent
(4) If $K$ is a Montesinos knot, then $m(K)\leq4$. 
\end{Proposition}

\vskip 3mm

\noindent{\bf Proof.} 
(1) It follows by Proposition \ref{one or two} (1). 

\noindent
(2) It follows by Proposition \ref{one or two} (2). 

\noindent
(3) Let $f=[(V,k)]$ be a morphism from $K$ to $T$ with $m_{\mathcal K}(f)=3$. Let $h:V\to{\mathbb S}^1\times{\mathbb D}^2$ be a homeomorphism with $m(h,k)=m_{\mathcal K}(f)=3$. Suppose that there are no critical points of $\pi\circ h|_k$. Then we have ${\rm braid}(K)\leq 3$. By Proposition \ref{braid and bridge} (1) ${\rm braid}(K)\leq2$ implies $m(K)\leq2$. Therefore ${\rm braid}(K)=3$. Suppose that there are some critical points of of $\pi\circ h|_k$. Then the local maximum and local minimum must appear alternatively along ${\mathbb S}^1$. Then we have a connected sum of some $2$-bridge knots. 
The converse is similar. 

\begin{figure}[htbp]
      \begin{center}
\scalebox{0.5}{\includegraphics*{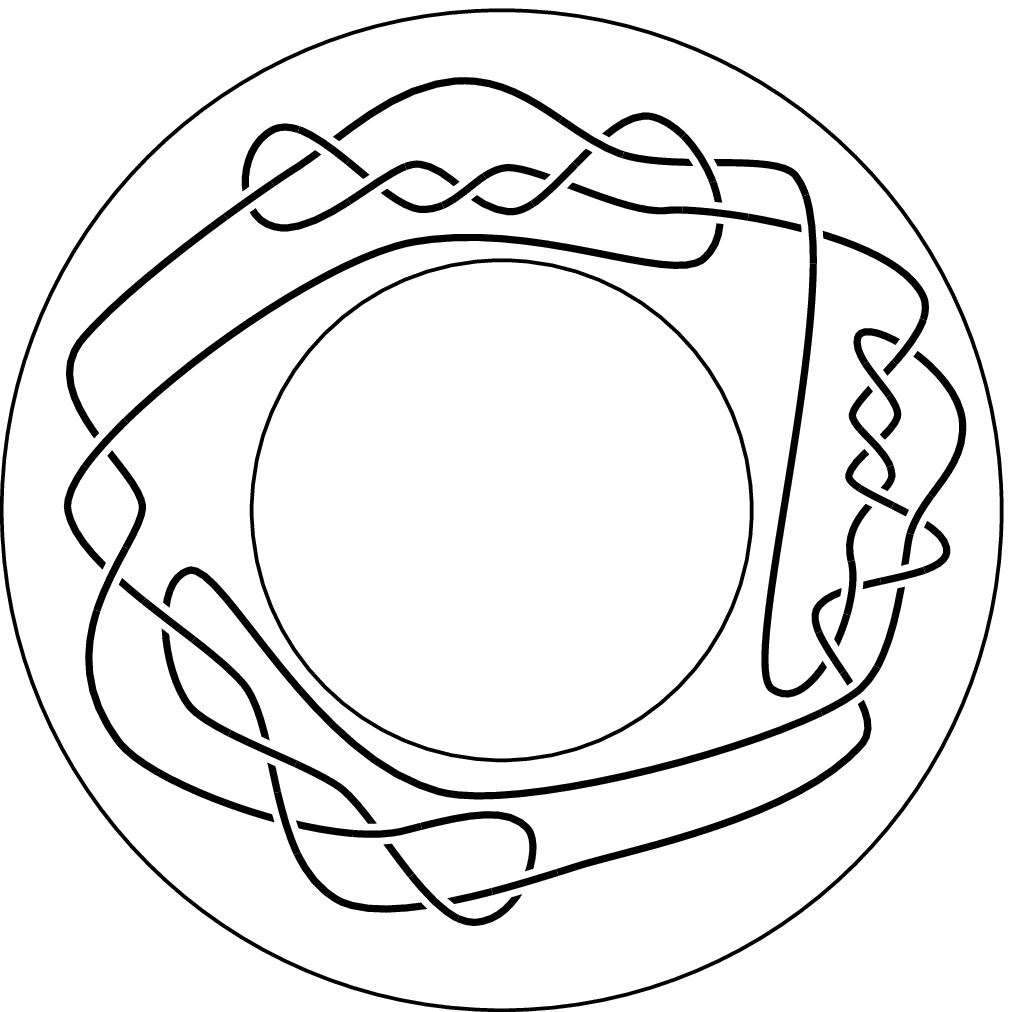}}
      \end{center}
   \caption{}
  \label{Montesinos}
\end{figure}

\noindent
(4) For a Montesinos knot $K$ it is easy to see $m(K)\leq4$. See for example Figure \ref{Montesinos}. 
$\Box$

\section*{Acknowledgments} The author is most grateful to Professors Kazuaki Kobayashi, Shin'ichi Suzuki and Ryuichi Ito. The author is also grateful to Professors Kazuhiro Kawamura, Kazufumi Eto, Ken-ichiroh Kawasaki, Youngsik Huh, Makoto Ozawa and Ryo Nikkuni for their helpful comments.

{\normalsize
\end{document}